\documentclass[final,leqno]{siamltex704}
\usepackage{graphicx}
\usepackage{epstopdf}
\usepackage{amsmath}
\usepackage{amsfonts}
\usepackage{multirow}
\usepackage{multirow}
\usepackage{caption}
\usepackage{rotating}
\usepackage{float}
\usepackage{subfigure}
\usepackage{hyperref}
\usepackage{amssymb}%
\usepackage[applemac]{inputenc}
\usepackage{fancybox, color}

\def\dom{\mathop{\rm Dom}}

\def\N{{\mathbb{N}}}

\def\R{{\mathbb{R}}}

\title{Stochastic Difference-of-Convex Algorithms for Solving nonconvex optimization problems}

\author{Le Thi Hoai An\thanks{%
Universit\'e de Lorraine, LGIPM, F-57000 Metz, France (\texttt{hoai-an.le-thi@univ-lorraine.fr})} \and 
Huynh Van Ngai\thanks{%
University of Quy Nhon, Viet Nam
(\texttt{ngaivn@yahoo.com})
}
\and Pham Dinh Tao\thanks{%
Laboratory of Mathematics, INSA-Rouen, University of Normandie, 76801 Saint-\'Etienne-du-Rouvray Cedex,
France (\texttt{pham@insa-rouen.fr})}
\and Luu Hoang Phuc Hau\thanks {%
Universit\'e de Lorraine, LGIPM, F-57000 Metz, France (\texttt{hoang-phuc-hau.luu@univ-lorraine.fr})}
}

\begin{document}

\date{}
\maketitle

\begin{abstract} The paper deals with stochastic difference-of-convex functions (DC)  programs, that is, optimization problems whose the cost function is a sum of a lower semicontinuous DC function and the expectation of a stochastic DC function with respect to a probability distribution. This class of nonsmooth and nonconvex stochastic optimization problems plays a central role in many practical applications. Although there are many contributions in the context of convex and/or smooth stochastic optimization, algorithms dealing with nonconvex and nonsmooth programs remain rare. In deterministic optimization literature, the DC Algorithm (DCA) is recognized to be one of the few algorithms to solve effectively nonconvex and nonsmooth optimization problems. The main purpose of this paper is to present some new stochastic DCAs for solving stochastic DC  programs. The convergence analysis of the proposed algorithms is carefully studied, and numerical experiments are conducted to justify the algorithms' behaviors.
\end{abstract}
\pagestyle{myheadings}
\thispagestyle{plain}
\markboth{Le Thi Hoai An,  Huynh Van Ngai, Pham Dinh Tao, and Luu Hoang Phuc Hau}{ \sc Stochastic DC algorithms }
\section{Introduction}

We consider the single stage stochastic optimization problems of the form
\begin{equation}\label{SOP}
\alpha=\inf\{f(x):=\varPhi(x)+r(x):\;\;\; x\in\R^n\},
\end{equation}  
where $\varPhi:\R^n\to\R$ is the expectation of a stochastic loss function with respect to the probability distribution $\mathbb{P}$ of the complete probability space $(\Omega,\Sigma_\Omega,\mathbb{P})$ 
\begin{equation}\label{int-func}
\varPhi(x):=\mathbb{E}_{s\simeq\mathbb{P}}\varphi(x,s)=\int_{\Omega}\varphi(x,s)d\mathbb{P},
\end{equation}
and $r:\R^n\to\R\cup\{+\infty\}$ is an extended real valued lower semicontinuous function. In general, the probability distribution $\mathbb{P}$ is unknown. A particular case, when $r$ is the indicator function of a closed convex set $C\subseteq \R^n,$ the problem reduces to minimizing an \textit{expected loss function} $\varPhi$ over a closed convex set,
$$\inf\{\varPhi(x)=\mathbb{E}_{s\simeq \mathbb{P}}\varphi(x,s): x\in C \}.$$
Stochastic optimization problems play a key role in many fields of applied science: Statistics, signal processing, finance, machine learning, and data science,... (see e.g., \cite{Bertsekas, BBL, EW-88, Han, Pflug, Pf-Rom, Sha-Den-Rus, Ste-Chris} and references given therein). Since the pioneering work by Robbin and Monro in 1951 \cite{RM} for solving stochastic programs with smooth and strongly convex data,  a huge number of publications related to their method for solving (\ref{SOP}) have been produced in both theory and application aspects. Generally, there are two principal approaches to stochastic optimization problems, along with some variants combining these two.

The first approach is approximating the stochastic loss function by a deterministic function in some appropriate stochastic ways to produce an approximation problem. The approximation problem is then optimized by deterministic (or stochastic, as well) optimization methods. Found solutions of the approximate problem are then regarded as approximate solutions of the original problem \cite{EN-13, HdM, NJLS, Sha-Den-Rus, Sha, Sha-89}. A popular approximation method is the Monte-Carlo sample average approximation described briefly as follows. Let $s_1,s_2,...,s_m$ be independent, identically distributed realizations obtained from the probability distribution $\mathbb{P},$ then the expected loss function is approximated by $\varPhi(x)\approx \varPhi_m(x):=\frac{1}{m}\sum_{i=1}^m \varphi(x,s_i),$
and the approximate optimization  problem is formulated as
\begin{equation}\label{MC approx prob}
\min\{\varPhi_m(x)+r(x):\quad x\in\R^n\}.
\end{equation}
In this approach, (\ref{MC approx prob}) is usually a large-sum problem since the number of samples $m$ would be very large, especially in the era of big data. Therefore, deterministic approaches to (\ref{MC approx prob}) would be prohibitively expensive; meanwhile, stochastic approaches would be efficient alternatives (e.g. \cite{Ber-MP11, AMNB, Mai, SRB-17}).

The second approach is iteratively constructing stochastic approximations (SA) for the expectation quantities (such as $\varPhi$ and $\nabla \varPhi$) of the original problem (\ref{SOP}), and then performing a solution-update step (e.g. see \cite{Borkar, Da-Drus, EW-88, GL-13, GLZ-16, XY-15} and references therein). An advantage of this approach is that the computational cost per iteration is cheap. Besides, solutions found by these algorithms are also the solutions (global, stationary, etc.) of the original problem (\ref{SOP}). However, a disadvantage is that the practical convergence rates are relatively slow since the methods' variance is large.

To our knowledge, thus far, most stochastic optimization methods in both approaches in the literature have mainly been developed for solving smooth and/or convex stochastic programs. In the nonsmooth and nonconvex setting, such algorithms remain rare. We list here the main approaches for solving nonconvex stochastic optimization problems. Most of the following works require the smoothness of the problem, in part or in full. The first approach is Stochastic (sub)gradient-based methods which compute stochastic gradients and perform a gradient-like update at each iteration (e.g. \cite{GL-13,Bertsekas00}). This approach is a natural extension from convex optimization to nonconvex $L$-smooth functions. Moreover, the proximal operator can be used in each update step, resulting in some variants called Stochastic proximal (sub)gradient-based methods (for example, see the article by Davis and Drusvyatskiy \cite{Da-Drus}). The second approach is Stochastic MM (Majorization-Minimization), which is a stochastic extension of deterministic MM (see works by Razaviyayn et al. \cite{Razaviyayn} and Mairal \cite{Mairal}). At each iteration, a sample surrogate function is constructed as an upper bound of the sample objective function. Then, the current sample surrogate will be averaged with all the past sample surrogates to create an average surrogate which is minimized to obtain an updated optimization variable. Liu et al. \cite{Liu20} take a step further as they proposed a stochastic MM scheme for compound stochastic programs, where the two outermost functions in the nested structure are assumed to be convex and isotone in order to preserve the convex property of the inner surrogates. The third approach is Stochastic SCA (Successive Convex Approximation) \cite{Scutari} which is similar to Stochastic MM, where the sequence of approximation functions are convex but not necessarily upper bounds of sample objective functions. The fourth approach is stochastic DCA that aims to handle stochastic DC programs - a very large class of nonsmooth and nonconvex optimization problems. Due to the challenges of this class of problems, this line of research begins with some special cases (e.g. large-sum, smooth) then gradually extends and elaborates. Le Thi et al. \cite{lethi17b} first developed stochastic DCA for large-sum problems of nonconvex but $L$-smooth functions with $\ell_{2,0}$ regularization, and later extended to a more general class of large-sum nonsmooth DC programs \cite{AMNB}. Clearly, these algorithms no longer work on the general setting of the form (\ref{SOP}). Liu et al. \cite{Liu20b} proposed a stochastic algorithm based on DCA to solve a special class of two-stage stochastic programs with linearly bi-parameterized quadratic recourse. Nitanda and Suzuki \cite{Nitanda} proposed a stochastic proximal DCA for differentiable DC programs and provided for the first time non-asymptotic convergence rate to find $\epsilon$-accurate solutions. Xu et al. \cite{Xu-al} further improved this work by designing stochastic proximal DCA schemes for a larger class of DC programs where the smoothness condition can be partially relaxed. 
It should be mentioned that in the above articles \cite{Razaviyayn, Mairal} dealing with stochastic MM  the authors also considered DC surrogates whose the second DC component is differentiable.

It is worth noting that, similar to deterministic optimization, most algorithms of the stochastic gradient-based, stochastic MM, stochastic SCA approaches (with usual choices of surrogates/approximation functions) can be seen as versions of stochastic DCA. Indeed, as indicated in \cite{LTP18}, even if the DC structure of the problem under consideration in existing approaches is hidden, the usual choices of the surrogate (reps. approximate) functions in MM (resp. SCA) methods result in DCA versions. Also, (proximal) (sub)gradient-based algorithms usually fall into the spectrum of DCA thanks to hidden DC structures of the problems at hand.
  
In this contribution, we are interested in \textit{Stochastic DC (SDC for short) } optimization problems (\ref{SOP}) with $\varphi(\cdot,s)$ given by: 
\begin{equation}\label{SDC func}
\varphi(x,s)=g(x,s)-h(x,s),\quad(x,s)\in\R^n\times\Omega,
\end{equation} 
where $g(\cdot,s)$ and $h(\cdot,s)$ are continuous convex functions, and $r$ is an extended real valued lower semicontinuous DC function. As $\varphi(\cdot,s)$ are DC, so is the expected loss $\varPhi$. Therefore, the original problem (\ref{SOP}) is naturally a DC program. However, the expected loss function is usually either unknown since the probability distribution $\mathbb{P}$ is unknown or too expensive to be exactly computed. This class of problems is very broad to cover almost all stochastic programs appearing in practice. To our knowledge, this is the first time in the literature such a general model is being considered. Indeed, we work with the general distribution $\mathbb{P}$, and allow both DC components of $\varphi$ to be nonsmooth. Moreover, the regularization term $r$ is also a nonsmooth DC function.
\vskip 0.1cm  
In the deterministic optimization literature, as DC optimization problems appear in many practical situations, DC programming plays a central role in nonconvex programming. The (deterministic) DCA was introduced in 1985 by Pham Dinh Tao \cite{Pham:chap:86} in the preliminary state and
extensively developed throughout various joint works by Le Thi Hoai An and Pham Dinh Tao (see  \cite{LTP18} and references therein)
to become now
classic and increasingly popular. Standard DCA solves DC programs of the form 
\begin{equation}\label{DCstandard}
f(x):=g(x)-h(x),\quad x\in\R^n,
\end{equation} 
where $g$ and $h$ (called DC components of $f$) are convex functions on $\R^n$. The main idea of DCA is quite simple: at each iteration $k$, DCA
approximates the second DC component $h(x)$ by its affine minorization  $h_k(x):=h(x^k)+\langle x-x^k,y^k\rangle$, 
with $%
y^{k}\in \partial h(x^{k})$, and minimizes the resulting convex function. 
\\
Nowadays, it is recognized that DCA is one of a few algorithms to solve effectively  nonconvex and nonsmooth programs, and there is a huge range of  applications of DCA in various fields of applied sciences. The DCA  was successfully applied to a lot of different optimization problems, and many nonconvex programs to  
which it gave almost always global solutions and was proved to be more robust
and more efficient than related standard methods, especially in the
large-scale setting. It is worth noticing that (see \cite{LTP18}) with
appropriate DC decompositions,   and suitably equivalent DC reformulations, DCA makes it possible to recover all (resp. most) standard methods in convex (resp. nonconvex) programming.
For instance, the readers are referred to   
\cite{An, LTSIAM2, ATao05, ATao-10, ANT, ANT-DCgeneral, ATMT, Pang-MOR,  PhLe1, PhLe2}, as well as   \cite{LTP18} for a survey on thirty years of developments of DC programming and DCA and references therein, and  very recent papers (e.g. \cite{APX, CPS, Liu20b, LPT-COA19, NPR, PT, QCLP}) for   nice properties of DC programming, DCA and their fruitful applications. 

In this paper, we develop stochastic DC algorithms for solving stochastic DC optimization problems of the form (\ref{SOP}) with $\varphi$ given by (\ref{SDC func}). 
Our proposed algorithms' features are very new compared to related works in the literature, which are highlighted as follows. Based on DCA, our main idea is to iteratively and randomly approximate the second DC component of the objective function as well as its subgradient while the first DC component is either approximated or left unchanged. 
We propose the following two variants of the Stochastic DC Algorithms (SDCA for short):
\begin{itemize}
\item SDCA with storage of past samples and subgradients;
\item SDCA with storage of past samples but updating subgradients. 
\end{itemize} 
For each variant, we develop two algorithms: the first algorithm iteratively and randomly approximates subgradients of the second DC component while keeping the first DC component unchanged; meanwhile, the second algorithm goes a step further as it also iteratively and randomly approximates the first DC component. In total, four stochastic algorithms are proposed.

The paper is organized as follows. In Section 2, we recall some basic notations and tools from Convex, Nonsmooth, and Variational Analysis, which will be used in the subsequent sections. Furthermore, we give a brief presentation on DC programming and  DCA. In Section 3, we present the  SDCA schemes and  their convergence results. Numerical experiments are presented in Section 4, while some concluding remarks and further research are discussed in the final section.
\section{Preliminaries}
\subsection{Tools from Convex and Variational Analysis }
Firstly we recall some notions from Convex Analysis and Nonsmooth Analysis, which will be needed thereafter (see, e.g., \cite{Borisbook1}, \cite{Roc}, \cite{R-W}).  In the sequel, the space $\R^n$ is equipped with the canonical inner product $\langle ,\rangle.$ Its dual space is identified with $\R^n$ itself. The open and closed balls with the center $ x\in \R^n$ and radius $\varepsilon >0$ are denoted, respectively, by $B(x,\varepsilon)$ and $B[x,\varepsilon],$ while the closed unit ball is denoted by $\mathbb{B}.$ 
 A function $f:\R^n\to\R\cup\{+\infty\}$ is called $\rho-$convex for some $\rho\ge 0,$  if for all $x,y\in\R^n,$ $\lambda\in[0,1]$ one has
$$f(\lambda x+(1-\lambda)y)\le \lambda f(x)+(1-\lambda) f(y)-\frac{\rho}{2}\lambda(1-\lambda)\|x-y\|^2.$$
The supremum of all $\rho\ge 0$ such that the above inequality holds is called the strong convexity modulus of $f,$ which is denoted by $\rho(f).$
\vskip 0.1cm
The conjugate of a convex function $f$ is denoted $f^*$ and is defined by
\begin{equation}\label{Conjugate}
f^{\ast }(y):=\sup \{\langle x,y\rangle -f(x):\,x\in \Bbb{R}^{n}\},\;\; y\in\R^n.
\end{equation}
The effective domain of $f$, denoted $\dom f$, is given by $\dom f:=\{x\in\R^n:\;\;f(x)<+\infty\}.$  The subdifferential of a convex function $f$ at $x\in\dom f$ is defined by
$\partial f(x)= \{x^*\in \R^n:\quad\langle x^*,y-x\rangle\le f(y)-f(x),\quad\forall y\in\R^n\}.$
We set $\partial f(x)=\emptyset$ if $x\notin\dom f.$ For a lower semicontinuous real extended valued function $f:\R^n\rightarrow\R\cup\{+\infty\}$, the \textit{Fr\'echet subdifferential} of $f$ at $x\in\dom f$ is defined by
$$\partial^Ff(x)=\left\{x^*\in\R^n:\;\; \liminf_{h\to 0}\frac{f(x+h)-f(x)-\langle x^*,h\rangle}{\|h\|}\geq 0\right\}.$$
For $x\notin\dom f,$ we set $\partial^F f(x)=\emptyset.$ The \textit{limiting subdifferential} of $f$ at $x\in\dom f$ is 
$$\partial f(x)=\{x^*\in\R^n:\;\; \exists (x_k,f(x_k))\to(x,f(x)),\; x_k^*\in\partial^Ff(x_k)),\;\;(x_k^*)\to x^*\},$$
and we put $\partial f(x) = \emptyset$ if $x\notin\dom f.$
If $f$ is locally Lipschitz at $x \in \mathbb{R}^n$, then the Clarke directional derivative $f^C(x,\cdot)$ at $x$ and the Clarke subdifferential $\partial^C f(x)$ are defined as \begin{align*}
    &f^C(x;d) = \lim_{(t,u) \to (0^+,x)} \dfrac{f(u+td)-f(u)}{t}\\
    \text{and } &\partial^Cf(x) = \{y \in \mathbb{R}^n: \langle y,d \rangle \leq f^C(x;d), \forall d \in \mathbb{R}^n\}.
\end{align*}

A point $x\in\R^n$ is called a \textit{Fr\'{e}chet (resp. limiting/ Clarke}) \textit{critical} point for the function $f,$  if $0\in\partial^Ff(x)$ (resp. $0\in \partial f(x)$/ $0 \in \partial^C f(x)$).
\vskip 0.1cm
When $f$ is a convex function, the Fr\'{e}chet, limiting, and Clarke subdifferential coincide with the subdifferential in the sense of Convex Analysis.

Let us recall the well-known sudifferential characterizations of the $\rho-$convexity (see e.g., \cite[Thm 5.1]{CSW-95}, \cite[Thm 8]{NP-08}).
\begin{theorem}\label{rho-conv Characs} Let $f:\R^n\to\R\cup\{+\infty\}$ be a lower semicontinuous function. For $\rho\ge 0$, the following  three statements are equivalent.
\begin{itemize}
\item[(i)] $f$ is a $\rho-$convex function.
\item[(ii)] For all $x,y\in\R^n$, $x^*\in \partial^Ff(x),$ one has
$\langle x^*,y-x\rangle\le f(y)-f(x)-\frac{\rho}{2}\|y-x\|^2.$
\item[(iii)] The sudifferential operator of $f,$ $\partial^Ff,$ is a $\rho-$monotone operator: for all $x,y\in\R^n,$ $x^*\in\partial^Ff(x),$ $y^*\in\partial^Ff(y),$
$\langle x^*-y^*,x-y\rangle\ge \rho\|x-y\|^2.$
\end{itemize}
\end{theorem}
 \subsection{A brief presentation on DC programming and DCA}
Let $\Gamma _{0}(\mathbb{R}^{n})$ denote the class of all lower semicontinuous proper real extended valued convex functions defined on $\mathbb{R}^{n}.$ The class
of DC functions is denoted by $DC(\mathbb{R}^{n})=\Gamma _{0}(\mathbb{R}^{n})-\Gamma _{0}(%
\mathbb{R}^{n})$, that is quite large to contain almost all real-life objective functions and is closed under all operations usually considered in
Optimization (see, e.g., \cite{ATao05}). We consider a standard DC program:
$$(\mathcal{P})\quad \alpha=\inf\{f(x):=g(x)-h(x):\;\; x\in\R^n\},$$
where $g,$ $h$ belong to $\Gamma_0(\R^n).$ Recall the natural convention $+\infty -(+\infty
)=+\infty $ and the fact  that if the optimal value $\alpha$ is finite, then $\dom g\subset \dom h.$
The dual problem of $(\mathcal{P})$ is defined by
$$(\mathcal{D})\quad \inf\{h^*(y)-g^*(y):\;\; y\in\R^n\},$$
where $g^{\ast }, h^{\ast}$ are the conjugate functions of $g,$ $h,$  respectively. Due to the duality result by Toland \cite{Toland} (see also \cite{PhLe1}), the optimal values of the primal dual problems coincide and there is the perfect symmetry between primal and dual programs
($\mathcal{P}$) and ($\mathcal{D}$): the dual program to ($\mathcal {D}$) is exactly ($\mathcal{P}$).
\vskip 0.1cm
A point $x^*\in\R^n$  is called  \textit{a DC critical} point of the DC problem  ($\mathcal{P}$) if $0\in\partial g(x^*)-\partial h(x^*),$ or  equivalently $\partial g(x^*)\cap\partial h(x^*) \neq \emptyset,$ while it is called \textit{a strongly DC critical point} of  ($\mathcal{P}$) if $\emptyset \neq \partial h(x^*) \subset \partial g(x^*).$ In the framework of DC programming, the terminology ``critical point" is referred to the notion of DC criticality. The notion of DC criticality is close to Clarke, Fr\'echet, limiting stationarity in the sense $\partial^F f(x^*) \subseteq \partial f(x^*) \subseteq \partial g(x^*) -\partial h(x^*)$ (whenever $h$ is continuous at $x^*$) and $\partial^C f(x^*) \subseteq \partial g(x^*)-\partial h(x^*)$, where equalities hold under technical assumptions. Furthermore, the commonly used \textit{directional stationarity} is equivalent to strong DC criticality, which is the strongest necessary condition for local DC optimality \cite{LTP18}.
\vskip 0.1cm
For a DC optimization problem  over a closed convex set constraint $C \neq \emptyset$, we can equivalently transform it into a standard DC program by using the indicator function of $C$ as
$\inf\{f(x)+\chi_C(x):\quad x\in \R^n\},$
where $\chi _{C}$ stands for the indicator function of $C,$ that is, $\chi _{C}(x)=0$ if $x\in C$, and $+\infty$ otherwise. For general DC programs with equality/ inequality constraints defined by DC functions, some penalty  techniques have been used to transform them to standard DC programs (see \cite{ANT-JOGO12, ANT-JOTA16}).  
\vskip 0.1cm
The DCA which is based on local optimality
and DC duality, consists in the construction of the two
sequences $\{x^{k}\}$ and $\{y^{k}\}$ (candidates for being primal and dual solutions, respectively) such that the
sequences of values of the primal and dual objective functions $\{g(x^{k})-h(x^{k})\}$, $\{h^{\ast }(y^{k})-g^{\ast }(y^{k})\}$
are decreasing, and their corresponding limits $%
x^{\infty }$ and $y^{\infty }$ satisfy local optimality conditions (see, e.g., \cite{LTSIAM2}, \cite{ATao05}, \cite{PhLe1}, \cite{PhLe2}). Briefly, the standard DCA is described as follows. Starting a given $x^0\in\dom g,$ and for $k=0,1,...,$ set
\begin{equation}\notag
\textrm{(DCA)}\quad \quad y^k\in \partial h(x^k);\quad x^{k+1}\in\partial g^*(y^k)=\textrm{argmin}\{g(x)-\langle y^k,x\rangle:\quad x\in\R^n\}.
\end{equation}
\section{ Stochastic DC Algorithms and convergence analysis}
Let $(\Omega,\Sigma_\Omega,\mathbb{P})$ be a probability space. Consider the stochastic DC program:
\begin{equation}\label{SDCP}
\alpha=\min\{f(x):=\varPhi(x)+r(x):\quad x\in\R^n\},
\end{equation}
where, $r:\R^n\to\R\cup\{+\infty\}$ is a lower semicontinuous DC function given by
\begin{equation}\label{r func}
r(x):=r_1(x)-r_2(x),\;\;x\in\R^n,
\end{equation} where $r_1,r_2:\R^n\to\R\cup\{+\infty\}$ are lower semicontinuous convex functions, and the expected loss function
\begin{equation}\label{Dc phi func}
\varPhi(x)=\mathbb{E}_s [\varphi(x,s)=g(x,s)-h(x,s)]=\int_\Omega[g(x,s)-h(x,s)]d\mathbb{P},
\end{equation}
with respect to continuous convex functions $g(\cdot,s)$, $h(\cdot,s),$ $s\in\Omega,$ defined on $\R^n.$
Throughout the paper, we assume that the expectations of $g(x,\cdot)$ and $h(x,\cdot)$ are finite for all $x\in\R^n,$ and denoted by
\begin{equation}\label{Ex convex func}
G(x):=\mathbb{E}_s[g(x,s)]=\int_{\Omega} g(x,s)d\mathbb{P},\quad H(x):=\mathbb{E}_s[h(x,s)]=\int_{\Omega} h(x,s)d\mathbb{P}.
\end{equation}
Then, $G,H$ are continuous convex functions on the whole space $\R^n,$ and therefore the objective function of the problem (\ref{SDCP}) admits a DC decompsition: $f=(G+r_1)-(H+r_2).$ In what follows we will make use of the following assumptions:
\begin{itemize}
\item[(A1)] For each $x\in\R^n,$ the function $h(x,\cdot)$ is bounded on $\Omega$ and  the functions $h(\cdot,s)$ are  equi-continuous at each $x\in \R^n,$ for $s\in\Omega,$ that is, for each $x\in\R^n,$ for any $\varepsilon>0,$ there exists $\delta>0$ such that
$|h(z,s)-h(x,s)|\le\varepsilon,\;\;\forall z\in B(x,\delta)$, and  $\forall s \in \Omega.$
\item[(A2)] For each $x\in\R^n,$ the functions $g(x,\cdot)$, $h(x,\cdot)$ are bounded on $\Omega$ and  the functions $g(\cdot,s)$, $h(\cdot,s)$ are equi-continuous at each $x\in \R^n,$ for $s\in\Omega.$
\end{itemize}
Recall that a critical point $x^*\in\R^n$ of problem (\ref{SDCP}) is characterized as follows 
\begin{equation}\label{Crit point}
0\in \partial (G+r_1)(x^*)-\partial (H+r_2)(x^*)=\partial G(x^*)+\partial r_1(x^*)-\partial H(x^*)-\partial r_2(x^*).
\end{equation}
Therefore, for $x^*\in\R^n$, the distance 
\begin{equation}\label{Dist to Stationary}
\begin{array}{ll}&d(0,\partial G(x^*)+\partial r_1(x^*)-\partial H(x^*)-\partial r_2(x^*))\\
&=\inf\{\|y\|:\;\;y\in \partial G(x^*)+\partial r_1(x^*)-\partial H(x^*)-\partial r_2(x^*)\},
\end{array}
\end{equation}
serves as a measure of ``proximity to criticality".
\vskip 0.1cm
Many practical optimization models in various fields of science and engineering can be formulated as stochastic DC optimization problems (\ref{SDCP}). Note that this class of programs (SDC) contains convex-composite and weakly convex optimization problems.  For the sake of illustration, let us give an example of the following robust real phase retrieval problem (see e.g., \cite{{Ca-Li-Sol}, Da-Drus, Du-Ru}):
\begin{equation}\label{Retrieval}
\min\{\mathbb{E}_{a,b}|\langle a,x\rangle^2-b|:\;\; x\in\R^n\},
\end{equation}
where, $a\in\R^n$, $b\in\R$ are independent random variables with given probability distributions. Usually, $a$ is a standard Gaussian random vector in $\R^n,$ and $b$ is defined by
$b=\langle a,\bar x\rangle^2+ \eta,$ with a noise $\eta.$ Obviously, the functions $\varphi(\cdot,a,b):=|\langle\cdot,a\rangle^2-b|$
are  DC, therefore, problem (\ref{Retrieval}) belongs to the class of stochastic DC programs. In particular, when $a,b$ are random variables of  the uniform distribution on a finite set of $m$ elements  with the  values respectively $\{a_1,...,a_m\}$ and $\{b_1,...,b_m\},$  the problem (\ref{Retrieval}) reduces to the following optimization problem
\begin{equation}\label{Model 1}
\min \left\{f(x)=\frac{1}{m}\sum_{i=1}^m|\langle a_i,x\rangle^2-b_i|:\quad x\in\R^n\right\}.
\end{equation}
In another perspective, the problem (\ref{Model 1}) can be regarded as an approximate problem of $(\ref{Retrieval}),$ when $a_i$ and $b_i$ ($i=1,...,m$) are realizations of $a$ and $b$, respectively. 
The latter problem can be reinterpreted as follows: Find $x\in\R^n$ such that $b_i \approx \langle a_i,x\rangle^2,\; i=1,2,...,m.$ 
By noticing the functions $\langle a_i,x\rangle^2-b_i$ are convex functions, for each $i=1,...,m,$ the function $|\langle a_i,x\rangle^2-b_i|$ admits a DC decomposition as follows,
$$|\langle a_i,x\rangle^2-b_i|=2\max\{\langle a_i,x\rangle^2-b_i,0 \}-(\langle a_i,x\rangle^2-b_i):= g_i(x)-h_i(x),\quad x\in\R^n,$$
where, for $i=1,...,m, ~g_i(x):=2\max\{\langle a_i,x\rangle^2-b_i,0 \};\;\;h_i(x):=\langle a_i,x\rangle^2-b_i,\; x\in\R^n.$
Then, a DC decomposition of the objective function is
$$f(x)=\frac{1}{m}\sum_{i=1}^m g_i(x)-\frac{1}{m}\sum_{i=1}^m h_i(x):=G(x)-H(x).$$ 
We are now presenting our proposed Stochastic DC Algorithms.
\subsection{Algorithms with storage of past samples and subgradients}
Firstly, we propose the following two algorithms in which at each iteration, the realized samples as well as subgradients from the past iterations are inherited. It is worth noticing that in these first algorithms, at each iteration, we have to compute only one subgradient of the function $h(\cdot,s)$ with respect to one current realization of $s$. 
We pick a sequence of positive reals $\{\alpha_k\}$ with $\sum_{k=0}^\infty\alpha_k=+\infty$. The first algorithm iteratively and randomly approximates the subgradient of the second DC component while keeping the first DC component unchanged.

\noindent\rule{12.7cm}{1.5pt}\\
\noindent \textbf{Algorithm 1:  Stochastic DC Algorithm 1 (SDCA1) }\\
\noindent\rule{12.7cm}{1pt}\\
\texttt{Initialization:} Initial data: $x^0\in\dom r_1,$ draw $s^0\overset{\textrm{iid}}{\simeq} \mathbb{P},$ and set $k=0.$ 
 \vskip 0.1cm
\texttt{Repeat:}
\begin{itemize}
 \item[1.] Compute $z^{k}\in \partial h(x^k,s^k);$ $u^k\in\partial r_2(x^k).$
 
 \item [2.] Set $y^k=\frac{1}{\sum_{i=0}^k\alpha_i}\sum_{i=0}^k\alpha_iz^i, \;\;w^k=\frac{1}{\sum_{i=0}^k\alpha_i}\sum_{i=0}^k\alpha_iu^i.$

\item [3.]  Compute a solution $x^{k+1}$ of the convex program
\begin{equation}\label{CP1}
\min \{G(x)+r_1(x)-\langle y^{k}+w^k,x\rangle :\ x\in \Bbb{R}^{n}\}.
\end{equation}
\item[4.] Set $k:=k+1$ and draw $s^k\overset{\textrm{iid}}{\simeq} \mathbb{P}.$
\end{itemize}

\texttt{Until} Stopping criterion.

\noindent\rule{12.7cm}{1.5pt}
\vskip 0.1cm

\bigskip

As observed, the first scheme implicitly supposes that the first DC component is in explicit form. To cope with the general situation where both DC components are unknown, we further introduce the stochastic DC algorithm 2, which differs from SDCA1 only by the replacement $G$ by its approximation $G_k$.

\noindent\rule{12.7cm}{1.5pt}\\
\noindent \textbf{Algorithm 2: Stochastic DC Algorithm 2 (SDCA2) }\\
\noindent\rule{12.7cm}{1.5pt}\\

Similar to Algorithm 1, where $G$ in step 3 of Algorithm 1 is replaced by $G_k$,
 $$G_k(x):=\frac{1}{\sum_{i=0}^k\alpha_i}\sum_{i=0}^k\alpha_ig(x,s^i),\;\;x\in\R^n.$$

\noindent\rule{12.7cm}{1.5pt}
\vskip 0.1cm
The convergence of Algorithms 1 and 2 is stated in the next theorem. Firstly, we need the following lemmas.
\begin{lemma}\label{Conver Diver series} For any increasing sequence of positive reals $\{\gamma_k\}$ with $\lim_{k\to\infty}\gamma_k=+\infty,$ one has $\sum_{k=1}^\infty(\gamma_k-\gamma_{k-1})/\gamma_k=+\infty.$
\end{lemma}
\vskip 0.1cm
\textit{Proof.}  By contradiction, we assume $\sum_{k=1}^\infty(\gamma_k-\gamma_{k-1})/\gamma_k<+\infty.$ Consequently, $\lim_{k\to\infty}\gamma_{k-1}/\gamma_k=1,$ and since $\lim_{t\to 1}(t-1)/\ln t =1,$ one has
$\lim_{k\to\infty}\frac{1-\gamma_{k-1}/\gamma_k}{\ln(\gamma_k/\gamma_{k-1})}=1.$
As $\sum_{k=1}^\infty\ln(\gamma_k/\gamma_{k-1})=\lim_{k\to\infty}(\ln\gamma_k-\ln \gamma_0)=+\infty,$
we derive the conclusion.\hfill{$\Box$}
\vskip 0.1cm 
The second is a variant of the strong law of large numbers. The proofs of this lemma and the next lemma will be given in Appendix. 
\begin{lemma}\label{STLN} Let $\{\alpha_k\}$ be a sequence of positive reals such that 
$$\frac{\sum_{i=0}^k\alpha_i^2}{\left(\sum_{i=0}^k\alpha_i\right)^2}\le \frac{N}{k^\gamma},\;\;\mbox{for all}\;k\in\N,\;\;\mbox{for some}\; N>0,\; \gamma>0.$$ Let $\{X_k\}$ be a sequence of independent and identically distributed (i.i.d. for short) random variables with $\mathbb{E}X_k=\mu.$ If either $\mathbb{E}X_k^4<+\infty$ and $\gamma>1/2$ or $\mathbb{E}X_k^2<+\infty$ and $\lim_{k,l\to\infty,l/k\to 1}\frac{\sum_{i=0}^l\alpha_i}{\sum_{i=0}^k\alpha_i}=1,$ then almost surely
$\frac{\sum_{i=0}^k\alpha_i X_i}{\sum_{i=0}^k\alpha_i}\rightarrow \mu\;\;\mbox{as}\;\; k\to\infty.$
\end{lemma}
\vskip 0.1cm
The next lemma is a weighted variant  of the uniform law of large numbers (\cite[ Lemma B.2]{EN-13}).
\begin{lemma}\label{ULLN} Let $X\subseteq \R^n$ be a compact set and let $(\Omega,\Sigma_\Omega,\mathbb{P})$ be a complete probability space. Let  $f: X\times \Omega \to\R$ be a function such that functions $f(\cdot,\omega)$ are uniformly bounded and equi-H\"{o}lder-continuous on $X$, that is, there are $M,L>0$ and $\gamma\in (0,1]$ such that
$|f(x,\omega)|\le M,\;\; |f(x,\omega)-f(y,\omega)|\le L\|x-y\|^\gamma,\;\forall x, y\in X,\; \omega \in\Omega.$ (It is noted that when $\gamma = 1,$ functions $f(\cdot,w)$ are called equi-Lipschitz.) Then there exists some constant $c>0$ such that for any sequence of positive reals $\{\alpha_k\}$ and any sequence of independent and identically distributed random variables with the probability distribution $\mathbb{P},$ $\{s^k\}$, one has
$$\mathbb{E}\max_{x\in X}\left|\frac{1}{\sum_{i=0}^k\alpha_i}\sum_{i=0}^k\alpha_i f(x,s^ i)-\mathbb{E}_sf(x,s)\right|\le \frac{c(1+\sqrt{\ln \beta_k})}{\beta_k},\;\;\mbox{for all}\;k\in\N_*,$$
where $\beta_k:=\frac{\sum_{i=0}^k\alpha_i}{\left(\sum_{i=0}^k\alpha_i^2\right)^{1/2}},\;\;k\in\N.$
\end{lemma}
\begin{theorem}\label{General Convergence} Let a sequence of positive reals $\{\alpha_k\}$ such that 
$\sum_{k=0}^\infty\alpha_k=+\infty,$
and for some $N,\gamma>0;$ for $\beta_k$ defined the same as in Lemma \ref{ULLN},
\begin{equation}\label{Cond Alpha seq}
 \frac{\sum_{i=0}^k\alpha_i^2}{\left(\sum_{i=0}^k\alpha_i\right)^2}\le \frac{N}{k^\gamma},\;\;\mbox{for all}\;\;k\in\N_*=\N\setminus\{0\}\quad\mbox{and}\quad\sum_{k=1}^\infty\frac{\alpha_k(1+\sqrt{\ln \beta_{k-1}})}{\beta_{k-1}\sum_{i=0}^k\alpha_i}<+\infty.
\end{equation} 
Assume further that either $\gamma>1/2$ or 
\begin{equation}\label{Quotient limit}
\lim_{k,l\to\infty,l/k\to 1}\frac{\sum_{i=0}^l\alpha_i}{\sum_{i=0}^k\alpha_i}=1,
\end{equation}
Suppose that (A1) holds for Algorithm 1; (A2) holds for Algorithm 2, and $\rho:=\inf_{s\in \Omega}\rho(h(\cdot,s))+\rho(r_2)>0.$ Let $\{x^k\}$ be a sequence generated by either Algorithm 1 or Algorithm 2. Suppose that the optimal value $\alpha$ of problem (\ref{SDCP}) is finite and with probability 1, $\limsup_{k\to\infty}\|x^k\|<\infty,$ and $\limsup_{k\to\infty}\|u^k\|<\infty.$  
 Then one has
\begin{itemize}
\item[(i)]  There exists a subsequence $\{x^{l_k}\}$ such that with probability 1, every its limit point is a critical point of problem (\ref{SDCP}).
\item[(ii)] (\textbf{The rate of the convergence to critical points}) Let $C$ be a compact set containing $\{x^k\}$. Suppose that for Algorithm 1, the derivative $\nabla r_2$ and the derivatives $\nabla h(\cdot,s), s \in \Omega$  are equi-Lipschitz with the same modulus $L/2$ on $C$; for Algorithm 2, the derivative $\nabla r_2$ and the derivatives $\nabla g(\cdot,s), \nabla h(\cdot,s) , s \in \Omega$  are equi-Lipschitz with the same modulus $L/2$ on $C$ . One has 
$d_k=O\left(1/\sqrt{\sum_{i=0}^k\alpha_i/A_i}\right)\;\;\mbox{ as}\quad k\rightarrow\infty,$
where, $A_k=\sum_{i=0}^{k}\alpha_i,\;k=0,1,\ldots$
 \end{itemize}
\end{theorem}
\vskip 0.1cm
\textit{Proof.}
\emph{Firstly, we prove Theorem \ref{General Convergence} for the Algorithm 1.}

$(i)$ Let the sequence $\{x^k\}$ be generated by Algorithm 1. Denote $A_k:=\sum_{i=0}^k\alpha_i,$ and define the functions $V_k:\R^n\to\R\cup\{+\infty\}$, $k\in\mathbb{N}$ by
$$V_k(x)=G(x)+r_1(x)-\frac{1}{A_k}\sum_{i=0}^{k}\alpha_i\langle z^i+u^i,x-x^ i\rangle-\frac{1}{A_k}\sum_{i=0}^k\alpha_i[h(x^ i,s^ i)+r_2(x^i)].$$
Denote by 
$\mathcal{F}_k=\sigma(s^0,...,s^{k-1}, x^0,...,x^{k},z^0,...,z^{k-1},u^0,...,u^{k-1}),\quad k\in\mathbb{N},$
the increasing $\sigma-$field generated by random variables $s^0,...,s^{k-1},x^0,...,x^{k},z^0,...,z^{k-1},$ and $u^0,...,u^{k-1}.$
As $x^{k+1}$ is a solution of the problem (\ref{CP1}), one has
\begin{equation}\label{estim 1}
\begin{array}{ll}
V_{k}(x^{k+1})&\le V_k(x^k)=G(x^k)+r_1(x^k)\\
&-\frac{1}{A_k}\sum_{i=0}^{k}\alpha_i\langle z^i+u^i,x^k-x^ i\rangle-\frac{1}{A_k}\sum_{i=0}^k\alpha_i[h(x^ i,s^ i)+r_2(x^i)].
\end{array}
\end{equation}
Next, due to the strong convexity of the function $h(\cdot,s)+r_2,$ $s\in\Omega$ with modulus at least $\rho,$ Theorem \ref{rho-conv Characs} implies
\begin{equation}\label{Estim 2}
\begin{array}{lll}
&V_{k-1}(x^k)-V_k(x^k)=\\
&\frac{\alpha_k}{A_k}[h(x^k,s^k)+r_2(x^k)]-\frac{\alpha_k}{A_kA_{k-1}}\left[\sum_{i=0}^{k-1}\alpha_i\left(\langle z^i+u^i,x^k-x^ i\rangle + h(x^i,s^ i)+r_2(x^i)\right)\right]\\
&\ge \frac{\alpha_k}{A_k}\left[h(x^k,s^k)-\frac{1}{A_{k-1}}\sum_{i=0}^{k-1}\alpha_ih(x^k,s^ i)\right]+\frac{\rho\alpha_k}{2A_kA_{k-1}}\sum_{i=0}^{k-1}\alpha_i\|x^k-x^ i\|^2.
 \end{array}
\end{equation}
From the preceding two inequalities, one has
\begin{equation}\label{Ge V-decreasing}
\begin{array}{ll}
V_{k}(x^{k+1})&\le V_{k-1}(x^{k})-(V_{k-1}(x^k)-V_k(x^k))\\
&\le V_{k-1}(x^{k})-\frac{\alpha_k}{A_k}\left[h(x^k,s^k)-\frac{1}{A_{k-1}}\sum_{i=0}^{k-1}\alpha_ih(x^k,s^ i)\right]\\
&-\frac{\rho\alpha_k}{2A_kA_{k-1}}\sum_{i=0}^{k-1}\alpha_i\|x^k-x^ i\|^2,
\end{array}
\end{equation}
and therefore, by taking conditional expectation both sides with respect to $\mathcal{F}_k$,
\begin{equation}\label{Ge Estim-Expec}
\mathbb{E}_{\mathcal{F}_k}V_{k}(x^{k+1})\le V_{k-1}(x^{k})-\xi_k-\frac{\rho\alpha_k}{2A_kA_{k-1}}\sum_{i=0}^{k-1}\alpha_i\|x^k-x^ i\|^2,
\end{equation}
where,
$\xi_k=\frac{\alpha_k}{A_k}\left[H(x^k)-\frac{1}{A_{k-1}}\sum_{i=0}^{k-1}\alpha_ih(x^k,s^ i)\right].$

By assumption, the sequence $\{x^k\}$ is bounded almost surely. Therefore, thanks to assumption (A1), we can apply Lemma \ref{ULLN} (by noting that $h(\cdot,s)$ are convex, so equi-continuity implies equi-Lipschitz), for some $c>0,$
\begin{equation}\label{Ge expec xi}
\mathbb{E}|\xi_k|\le\frac{c\alpha_k(1+\sqrt{\ln \beta_{k-1}})}{\beta_{k-1}A_k},\;\;\mbox{for all}\; k= 1,2,...
\end{equation}
Thus by the second relation of (\ref{Cond Alpha seq}), \begin{equation}\label{Bounded expec sum}
 \sum_{k=0}^\infty\mathbb{E}|\xi_k|<+\infty.
 \end{equation}
 \vskip 0.2cm
Since the sequence $\{x^k\}$ is bounded almost surely, we can assume it is contained in a compact set $C\subseteq\R^n.$ Pick $R>0$ such that  $C\subseteq B_R:=\{x\in\R^n:\;\|x\|<R\}.$ 
 Denote by
$\sigma_k(x):=\frac{1}{A_{k}} \sum_{i=0}^{k}\alpha_ih(x,s^ i).$
Let $\mathbb{Q}^n\subseteq \R^n$ be the set of points with all rational coordinates. Then $\mathbb{Q}^n\cap \overline{B}_R$ is a countable set which is dense in $\overline{B}_R$. Denote
$\mathbb{Q}^n\cap \overline{B}_R:=\{z^1,z^2,...,z^k,...\},$ then for each $l=1,2,...,$
$\mathbb{P}(\{\sigma_k(z^l)\nrightarrow H(z^l)\})=0.$ Hence, denoting by $S$ the event 
$S:=\bigcap_{l=1}^\infty\{\sigma_k(z^l)\rightarrow H(z^l)\},$ one has
$$1\ge\mathbb{P}(S):=\mathbb{P}\left(\bigcap_{l=1}^\infty\{\sigma_k(z^l)\rightarrow H(z^l)\}\right)\ge 1-\sum_{l=1}^\infty\mathbb{P}\left(\{\sigma_k(z^l)\nrightarrow H(z^l)\}\right)=1,$$
so $\mathbb{P}(S)=1.$
\vskip 0.2cm 
 \noindent\textit{Claim 1. For all $x\in\overline{B}_R,$ for all $\{s^k\}_{k\in \mathbb{N}}\in S,$ $\lim_{k\to\infty}\sigma_k(x)=H(x).$}
 \vskip 0.2cm 
 When $\gamma>1/2$ or (\ref{Quotient limit}) is verified, then the strong law of large numbers with weighted averages stated in Lemma \ref{STLN} also holds for the sequence $\{h(x,s^k)\}$ for each $x\in \R^n.$ 
 Let $\{s^k\}_{k\in\N}\in S$ be given. Let $x\in \overline{B}_R,$ and let $\{\varepsilon_l\}_{l\in\N}$ be a sequence of positive reals converging to $0.$ Then there is a subsequence  
 $\{z^{k_l}\}\subseteq\mathbb{Q}^n\cap\overline{B}_R$ such that $|H(x)-H(z^{k_l})|<\varepsilon_l$ as well as $|h(x,s)-h(z^{k_l},s)|<\varepsilon_l$ for all $l\in\N,$ all $s\in \Omega.$ For $\{s^k\}_{k\in\N}\in S,$ and for $l\in \N,$ since $\sigma_k(z^{k_l})\rightarrow H(z^{k_l}),$ there is an index $K_l$ such that $|\sigma_k(z^{k_l})- H(z^{k_l})|<\varepsilon_l$ for all $k\ge K_l.$ Hence, for all $k \ge K_l$,
 $$|\sigma_k(x)- H(x)|\le |\sigma_k(z^{k_l})- H(z^{k_l})|+|\sigma_k(x)-\sigma_k(z^{k_l})|+|H(x)-H(z^{k_l})|<3\varepsilon_l,$$
 which shows that 
 \begin{equation}\label{Almost surely} 
 \lim_{k\to\infty}\sigma_k(x)=\lim_{k\to\infty}\frac{1}{A_{k}} \sum_{i=0}^{k}\alpha_ih(x,s^ i)= H(x),\quad\forall \{s^k\}_{k\in\N}\in S,\;\forall x\in \overline{B}_R,
 \end{equation}
 and Claim 1 is proved.
 \vskip 0.2cm
 \noindent\textit{Claim 2. Almost surely, the sequence $\{V_k(x^{k+1})\}$ is bounded below.}
 One has
 $$
 \begin{array}{ll}V_k(x^{k+1})&=G(x^{k+1})+r_1(x^{k+1})-\frac{1}{A_k}\sum_{i=0}^{k}\alpha_i[\langle z^i,x^{k+1}-x^ i\rangle+h(x^i,s^i)]\\
 &-\frac{1}{A_k}\sum_{i=0}^k\alpha_i[\langle u^i,x^{k+1}-x^ i\rangle+r_2(x^i)]\\
 &\ge G(x^{k+1})+r_1(x^{k+1})- \frac{1}{A_k}\sum_{i=0}^{k}\alpha_ih(x^{k+1},s^i)-r_2(x^{k+1}),
 \end{array}
 $$
 which implies
 \begin{equation}\label{Lower estimate-V}
 V_k(x^{k+1})\ge f(x^{k+1})-\underbrace{\left[\frac{1}{A_k}\sum_{i=0}^{k}\alpha_ih(x^{k+1},s^i)-H(x^{k+1})\right]}_{(Q)}.
 \end{equation}

Since $\{x^k\}$ is bounded, from Claim 1 and assumption (A1) we derive $(Q)$ is bounded. This in conjunction with $\{f(x^{k+1})\}_{k\in\N}$ being bounded below by $\alpha$ imply - almost surely- the sequence $\{V_k(x^{k+1})\}$ is bounded below by some constant $\bar{\alpha}.$
 \vskip 0.2cm
  By considering the sequence $\{V_k(x^{k+1})-\bar{\alpha}\}_{k\in\N}$ instead of $\{V_k(x^{k+1})\}_{k\in\N},$ invoking (\ref{Ge Estim-Expec}) and (\ref{Bounded expec sum}), thanks to the convergence theorem for nonnegative almost supermartigales (\cite[Thm 1]{Rob-Sie}, see also, \cite{Dembo}), one derives that almost surely the sequence $\{V_k(x^{k+1})\}$ converges and 
\begin{equation}\label{General squared series}
\sum_{k=1}^\infty\frac{\alpha_k}{A_kA_{k-1}}\sum_{i=0}^{k-1}\alpha_i\|x^k-x^ i\|^2<+\infty.
\end{equation} 
In view of Lemma \ref{Conver Diver series}, $\sum_{k=1}^\infty\alpha_k/A_k=+\infty;$ consequently, there is a subsequence of 
$\delta_k:=\frac{1}{A_{k-1}}\sum_{i=0}^{k-1}\alpha_i\|x^k-x^ i\|^2$ converging almost surely to $0.$ Next, by picking a subsequence and relabeling if necessary, without loss of generality, we can assume that $\delta_k\to 0.$  
By the convexity of the function $(\cdot)^2,$ it follows from Jensen inequality that
$\left(\frac{1}{A_{k-1}}\sum_{i=0}^{k-1}\alpha_i\|x^{_k}-x^ i\|\right)^2\le \frac{1}{A_{k-1}}\sum_{i=0}^{k-1}\alpha_i\|x^{k}-x^ i\|^2,$
which implies almost surely $\lim_{k}\frac{1}{A_{k-1}}\sum_{i=0}^{k-1}\alpha_i\|x^{k}-x^ i\|=0.$
\vskip 0.2cm
 Let $\{s^k\}_{k\in\N}\in S$ be given such that $\lim_{k}\frac{1}{A_{k-1}}\sum_{i=0}^{k-1}\alpha_i\|x^{k}-x^ i\|=0.$ Suppose $x^*\in\R^n$ is a limit point of the sequence $\{x^{k}\},$ with respect to $\{s^k\}_{k\in\N}\in S,$ say, there is a subsequence $\{x^{l_k}\}$ converging to $x^*.$ Then, from the above relation,
$\lim_{k}\frac{1}{A_{l_k-1}}\sum_{i=0}^{l_k-1}\alpha_ix^ i=\lim_{k}x^{l_k}=x^*.$
By the equi-Lipschitz property of $h(\cdot,s)$
with a Lipschitz constant $L$ on $\overline{B}_R$, $\|z^i\|\le L,$ $(i=0,1,...)$ and $H$ is also Lipschitz on $\overline{B}_R$ with the same constant $L.$
Therefore, one has
$\frac{1}{A_{l_k-1}}\left|\sum_{i=0}^{l_k-1}\alpha_i\langle z^i,x^{l_k}-x^ i\rangle\right|\le\frac{L}{A_{l_k-1}}\sum_{i=0}^{l_k-1}\alpha_i\|x^{l_k}-x^ i\|,$
and
\begin{align*}
\left|\frac{1}{A_{l_k-1}}\sum_{i=0}^{l_k-1}\alpha_ih(x^ i,s^ i)-H(x^{l_k})\right|&\le \left|\frac{1}{A_{l_k-1}}\sum_{i=0}^{l_k-1}\alpha_ih(x^*,s^ i)-H(x^*)\right|\\&+\frac{L}{A_{l_k-1}}\sum_{i=0}^{l_k-1}\alpha_i\|x^{l_k}-x^ i\| + 2L \Vert x^{l_k}-x^{*} \Vert .
\end{align*}
As $\lim_{k\to\infty}\|x^{l_k}-x^*\|=\lim_{k\to\infty}\frac{L}{A_{l_k-1}}\sum_{i=0}^{l_k-1} \alpha_i\|x^{l_k}-x^ i\|=0,$ and by Claim 1,  
$\lim_{k\to\infty}\left[\frac{1}{A_{l_k-1}}\sum_{i=0}^{l_k-1}\alpha_ih(x^*,s^ i)-H(x^*)\right]=0,$ 
one has
\begin{equation}\label{Ge Average h}
\lim_{k\to\infty}\frac{1}{A_{l_k-1}}\sum_{i=0}^{l_k-1}\alpha_i\langle z^i,x^{l_k}-x^ i\rangle=0;\;\;
\lim_{k\to\infty}\left[\frac{1}{A_{l_k-1}}\sum_{i=0}^{l_k-1}\alpha_ih(x^ i,s^ i)- H(x^{l_k})\right]=0.
\end{equation}
 By passing to a subsequence if necessary,  assume that 
$$\lim_ky^{l_{k}-1}=\lim_k\frac{1}{A_{l_k-1}}\sum_{i=0}^{l_k-1}\alpha_iz^i=y^*,\;\mbox{and}\;\lim_kw^{l_{k}-1}=\lim_k\frac{1}{A_{l_k-1}}\sum_{i=0}^{l_k-1}\alpha_iu^i=w^*.$$ 
Since $y^{l_{k}-1}+w^{l_{k}-1}\in\partial G(x^{l_k})+\partial r_1(x^{l_k}),$ passing to the limit, one obtains $y^*+w^*\in\partial G(x^*)+\partial r_1(x^*).$
Next, $w^{l_{k}-1}\in \frac{1}{A_{l_k-1}}\sum_{i=0}^{l_k-1}\alpha_i\partial r_2(x^ i)$ and $\Vert u^k \Vert \leq M$ for all $k \in \mathbb{N},$ the Jensen inequality implies
$$\begin{array}{ll}
\langle w^{l_k-1},x-x^{l_k}\rangle &=\frac{1}{A_{l_k-1}}\sum_{i=0}^{l_k-1}\alpha_i\langle u^i,x-x^ i\rangle+\frac{1}{A_{l_k-1}}\sum_{i=0}^{l_k-1}\alpha_i\langle u^i,x^ i-x^{l_k}\rangle\\
&\le r_2(x)-\frac{1}{A_{l_k-1}}\sum_{i=0}^{l_k-1}\alpha_ir_2(x^ i)+\frac{1}{A_{l_k-1}}\sum_{i=0}^{l_k-1}\alpha_i\|u^i\|\|x^ i-x^{l_k}\|\\
&\le r_2(x)-r_2\left(\frac{1}{A_{l_k-1}}\sum_{i=0}^{l_k-1}\alpha_ix^ i\right)+\frac{M}{A_{l_k-1}}\sum_{i=0}^{l_k-1}\alpha_i\|x^ i-x^{l_k}\|,
\end{array}$$
for all $x \in \mathbb{R}^n$. Passing to the limit as $k\to\infty,$ by $\frac{1}{A_{l_k-1}}\sum_{i=0}^{l_k-1}\alpha_ix^ i\rightarrow x^*$ as well as $\frac{M}{A_{l_k-1}}\sum_{i=0}^{l_k-1}\alpha_i\|x^ i-x^{l_k}\|\rightarrow 0,$ one derives $w^*\in \partial r_2(x^*).$ On the other hand, $y^{l_{k}-1}\in \frac{1}{A_{l_k-1}}\sum_{i=0}^{l_k-1}\alpha_i\partial h(x^ i,s^ i),$ then
\begin{equation}\label{end subdiff}
\begin{array}{ll}\langle y^{l_k-1},x-x^{l_k}\rangle &=\frac{1}{A_{l_k-1}}\sum_{i=0}^{l_k-1}\alpha_i\langle z^i,x-x^ i\rangle+\frac{1}{A_{l_k-1}}\sum_{i=0}^{l_k-1}\alpha_i\langle z^i,x^ i-x^{l_k}\rangle\\
&\le\frac{1}{A_{l_k-1}} \sum_{i=0}^{l_k-1}\alpha_ih(x,s^ i)-\frac{1}{A_{l_k-1}}\sum_{i=0}^{l_k-1}\alpha_ih(x^ i,s^ i)\\
&+\frac{L}{A_{l_k-1}}\sum_{i=0}^{l_k-1}\alpha_i\|x^ i-x^{l_k}\|\quad\forall x\in B_R .
\end{array}
\end{equation}
Noticing that $\frac{1}{A_{l_k-1}}\sum_{i=0}^{l_k-1}\alpha_i\|x^ i-x^{l_k}\|\rightarrow 0,$ and by (\ref{Ge Average h}), $\frac{1}{A_{l_k-1}} \sum_{i=0}^{l_k-1}\alpha_ih(x^ i,s^ i)\rightarrow H(x^*).$ Now, by letting $k\to\infty$ in the relation (\ref{end subdiff}), one obtains $y^*\in\partial H(x^*).$ Consequently $y^*+v^*\in \left(\partial G(x^*)+\partial r_1(x^*) \right)\cap \left( \partial H(x^*)+\partial r_2(x^*) \right)$. Thus,  $x^*$ is a critical point of problem (\ref{SDCP}).
\vskip 0.1cm
$(ii)$ Assume the derivative $\nabla r_2$ and the derivatives $\nabla h(\cdot,s),$ $s\in \Omega,$ are equi-Lipschitz with the same modulus $L/2$ on a compact set containing $\{x^k\}.$ Then, $H$ is differentiable and $\nabla H(x)=\mathbb{E}_s\nabla h(x,s),$ for all $x\in\R^n.$ For $k\in\N\setminus\{0\},$
\begin{equation}\label{Deriv-Lip}
\begin{array}{ll}
y^{k-1}+w^{k-1}&= \frac{1}{A_{k-1}}\sum_{i=0}^{k-1}\alpha_i[\nabla h(x^ i,s^ i)+\nabla r_2(x^i)]\\
&\in \frac{1}{A_{k-1}}\sum_{i=0}^{k-1}\alpha_i\nabla h(x^k,s^ i)+\nabla r_2(x^k)+\frac{L}{A_{k-1}}\sum_{i=0}^{k-1}\alpha_i\|x^k-x^ i\|\mathbb{B}\\
&\subseteq \nabla H(x^k) +\nabla r_2(x^k)+\left(\eta_k+L\mu_k\right)\mathbb{B}, 
\end{array}\end{equation}
where,
$\mu_k=\frac{1}{A_{k-1}}\sum_{i=0}^{k-1}\alpha_i\|x^k-x^ i\|$ and $ \eta_k=\left\|\frac{1}{A_{k-1}}\sum_{i=0}^{k-1}\alpha_i \nabla h(x^k,s^ i)-\nabla H(x^k)\right\|.$ As $y^{k-1}+w^{k-1}\in\partial G(x^k)+\partial r_1(x^k),$ we derive
\begin{equation}\label{dist-Critical}
d(0,\partial G(x^k)+\partial r_1(x^k)-\nabla H(x^k)-\nabla r_2(x^k))\le \eta_k+L\mu_k. 
\end{equation}
By using Lemma \ref{ULLN} for each vector component, there is, say, the same constant $c>0$ above,  such that for all $k\in\N_*,$ 
\begin{equation}\label{nabla-estim}
\mathbb{E}\eta_k=\mathbb{E}\left\|\frac{1}{A_{k-1}}\sum_{i=0}^{k-1}\alpha_i\nabla h(x^k,s^ i)-\nabla H(x^k)\right\|\le \frac{c(1+\sqrt{\ln \beta_{k-1}})}{\beta_{k-1}}.
\end{equation}
Inequalities (\ref{dist-Critical}) and (\ref{nabla-estim}) yield
\begin{equation}\label{Inq-Dist}
\tau_k:=\mathbb{E}d(0,\partial G(x^k)+\partial r_1(x^k)-\nabla H(x^k)-\nabla r_2(x^k))\le L \mathbb{E}\mu_k+ \frac{c(1+\sqrt{\ln \beta_{k-1}})}{\beta_{k-1}}.
\end{equation}
Consequently, by using the Cauchy inequality, $(a+b)^2\le 2(a^2+b^2)$ for $a,b\ge 0,$ 
$$(\mathbb{E}\mu_k)^2\ge \tau_k^2/(2L^2)- c^2(1+\sqrt{\ln \beta_{k-1}})^2/(L^2\beta_{k-1}^2),$$
and as $\mu^2_k\le\delta_k:=\frac{1}{A_{k-1}}\sum_{i=0}^{k-1}\alpha_i\|x^k-x^ i\|^2,$ and by using the relation $(\mathbb{E}\mu_k)^2\le \mathbb{E}\mu_k^2,$ from (\ref{Ge expec xi}) and (\ref{Ge V-decreasing}), one obtains 
$$\begin{array}{ll}\mathbb{E}V_{k}(x^{k+1})- \mathbb{E}V_{k-1}(x^k)&\le 
-\frac{\rho\alpha_k\mathbb{E}\delta_k}{2A_k}+\mathbb{E}|\xi_k|\le -\frac{\rho\alpha_k\mathbb{E}\mu^2_k}{2A_k}+\frac{c\alpha_k(1+\sqrt{\ln \beta_{k-1}})}{\beta_{k-1}A_k}\\
&\le -\frac{\rho\alpha_k\tau_k^2}{4L^2A_k}+ \frac{\rho c^2\alpha_k(1+\sqrt{\ln \beta_{k-1}})^2}{2L^2\beta_{k-1}^2A_k}+\frac{c\alpha_k(1+\sqrt{\ln \beta_{k-1}})}{\beta_{k-1}A_k}.
\end{array}$$ 
By adding these inequalities with $k=1,2,...,$ one derives, for $k\in\N_*,$
$$
\sum_{i=1}^k\frac{\rho\alpha_i\tau_i^2}{4L^2A_i}\le \mathbb{E}V_0(x^1)-\mathbb{E}V_{k}(x^{k+1})+\sum_{i=1}^k\left[\frac{\rho c^2\alpha_i(1+\sqrt{\ln \beta_{i-1}})^2}{2L^2 \beta_{i-1}^2A_i}+\frac{c\alpha_i(1+\sqrt{\ln  \beta_{i-1}})}{ \beta_{i-1}A_i}\right].
$$
Since $d_k\le\tau_i,$ for $i=1,...,k,$ and $\{V_k(x^{k+1})\}$ is bounded a.s., say $V_k(x^{k+1})\ge V^*$ for some $V^*\in\R,$ for all $k\in \N,$ the inequality above implies
$$\frac{\rho}{4L^2}d^2_k\le\frac{\mathbb{E}V_0(x^1)-V^*+\sum_{i=1}^k\left[\frac{\rho c^2\alpha_i(1+\sqrt{\ln  \beta_{i-1}})^2}{2L^2 \beta_{i-1}^2A_i}+\frac{c\alpha_i(1+\sqrt{\ln  \beta_{i-1}})}{ \beta_{i-1}A_i}\right]}{\sum_{i=1}^k\frac{\alpha_i}{A_i}}=O\left(\frac{1}{\sum_{i=1}^k\frac{\alpha_i}{A_i}}\right). $$ 

\emph{Secondly, we prove Theorem \ref{General Convergence} for the Algorithm 2.}

Consider the sequence $\{x^k\}$ generated by Algorithm 2. The proof is similar to the Algorithm 1's proof, so we sketch it. 

$(i)$ Consider the function $V_k:\R^n\to\R\cup\{+\infty\}$, $k\in\mathbb{N}$ defined as
$$V_k(x):=\frac{1}{A_k}\sum_{i=0}^k\alpha_ig(x,s^ i)+r_1(x)-\frac{1}{A_k}\sum_{i=0}^{k}\alpha_i\langle z^i+u^i,x-x^ i\rangle-\frac{1}{A_k}\sum_{i=0}^k\alpha_i[h(x^ i,s^ i)+r_2(x^i)].$$

By using the same arguments of the proof of Algorithm 1, we arrive at
\begin{align*}
V_{k}(x^{k+1})&\le V_{k-1}(x^{k})-\frac{\rho \delta_k\alpha_k}{2A_k}\\
&+\frac{\alpha_k}{A_k}\left[g(x^k,s^k)-h(x^k,s^k)-\frac{1}{A_{k-1}}\sum_{i=0}^{k-1}\alpha_i\left(g(x^k,s^ i)-h(x^k,s^ i)\right)\right],
\end{align*}
where $\delta_k:=\frac{1}{A_{k-1}}\sum_{i=0}^{k-1}\alpha_i\|x^{k}-x^ i\|^2.$
Taking expectation with respect to $\mathcal{F}_k,$ one obtains
$\mathbb{E}_{\mathcal{F}_k}V_{k}(x^{k+1})\le V_{k-1}(x^{k})-\frac{\rho \delta_k\alpha_k}{2A_k}+\xi_k,$
where
$$\xi_k:=\frac{\alpha_k}{A_k}\left[G(x^k)-H(x^k)-\frac{1}{A_{k-1}}\sum_{i=0}^{k-1}\alpha_i\left(g(x,s^ i)-h(x,s^ i)\right)\right];$$
The rest of the proof is completely similar to the preceding one, here instead of $h$, we make use of assumption (A2) to obtain the equi-Lipschitz property of the functions $g(\cdot,s)$ and $g(\cdot,s)-h(\cdot,s),$ as well. Thus, by the supermartigale theorem, we also derive the existence of a subsequence of $\{\delta_k\}$ converging to $0.$ As in the previous case one has $y^*+w^*\in\partial H(x^*)+r_2(x^*).$ To show $y^*+w^*\in\partial G(x^*)+\partial r_1(x^*),$ noticing that in this case,
$y^{l_k-1}+w^{l_k-1}\in\frac{1}{A_{l_k-1}}\sum_{i=0}^{l_k-1}\alpha_i\partial g(x^{l_k},s^i)+\partial r_1(x^{l_k}),$ one has
$$\langle y^{l_k-1}+w^{l_k-1},x-x^{l_k}\rangle\le \frac{1}{A_{l_k-1}}\sum_{i=0}^{l_k-1}\alpha_ig(x,s^i)+r_1(x)-\frac{1}{A_{l_k-1}}\sum_{i=0}^{l_k-1}\alpha_ig(x^{l_k},s^i)-r_1(x^{l_k}),$$
for all $x\in\R^n.$ By (A2), Claim 1 also holds for $g(\cdot,s),$ that is, almost surely, $\frac{1}{A_{l_k-1}}\sum_{i=0}^{l_k-1}\alpha_ig(\cdot,s^i)\to G(\cdot)$ on every bounded set in $\R^n.$ Hence by letting $k\to\infty$ in the previous inequality, we obtain $y^*+w^*\in\partial G(x^*)+\partial r_1(x^*).$

$(ii)$ 
Following the arguments of part $(ii)$ of Algorithm 1, we arrive at (\ref{dist-Critical}) with $G$ being replaced by $G_k$. Note that $\nabla g(\cdot,s), s \in \Omega$ are equi-Lipschitz, arguments applied to $h(\cdot,s)$ are applicable to $g(\cdot,s)$ to obtain the desired convergence rate.

 \hfill {$\Box$}
\vskip 0.2cm
  \noindent\textit{Remark 1.} $(i)$ Observe from the proof that to obtain the convergence rate in part $(ii)$ of the theorem, it needs only the convergence of the series $\sum_{k=1}^\infty\frac{\alpha_k(1+\sqrt{\ln \beta_{k-1}})}{\beta_{k-1}\sum_{i=0}^k\alpha_i}.$
  \vskip0.1cm
  $(ii)$ An usual way to choose the sequence $\{\alpha_k\}_{k\in\N}$ is to take $\alpha_k=1,$ for all $k\in\N.$ For this sequence $\alpha_k$, the convergence rate of  Algorithms 1 and 2 given in Theorem \ref{General Convergence} (ii) is 
  $d_k=O\left(1/\sqrt{\sum_{i=0}^k\frac{1}{i+1}}\right)=O(1/\sqrt{\ln k}).$
  \vskip 0.1cm
  $(iii)$ Obviously, the sequence $\alpha_k:=k^\alpha,\;k\in\N_*$ with $\alpha\ge-1/2$ verifies all three conditions (\ref{Cond Alpha seq}) and (\ref{Quotient limit}) of Theorem \ref{General Convergence}. Let us take another example of the sequence $\{\alpha_k\}$ which produces an asymptotic convergence rate better than the rate obtained by the sequence $\{k^\alpha\}$ ($\alpha \geq -1/2).$ For $a>1$ and $\alpha\in (0,1/2),$ let $\{\alpha_k\}$ be a sequence of positive reals such that
  $0<\lim_{k\to\infty}\frac{A_k}{a^{k^\alpha}}<+\infty,$
 where $A_k=\sum_{i=0}^k \alpha_i.$ For example, one can take the sequence $\alpha_k=a^{(k+1)^\alpha}-a^{k^\alpha},\;k\in\N_*.$
 One sees obviously that $\alpha_k=O(a^{k^\alpha}k^{\alpha-1}),$ therefore
  $\label{Rate estimmate}
  \alpha_k/A_k=O\left(k^{\alpha-1}\right),\;\;\mbox{and}\;\; \sum_{i=0}^k\alpha_i/A_i=O(k^\alpha).$
  One has the following estimate
  $$\begin{array}{ll}
  \sum_{i=1}^k\alpha_i^2\le C\sum_{i=1}^{k}a^{2i^\alpha}i^{2(\alpha-1)}\le Ca^{2k^\alpha}\sum_{i=1}^ki^{2(\alpha-1)}\le C_1a^{2k^\alpha}k^{2\alpha-1},
  \end{array}$$
  for some $C_1>C>0.$ It implies the first relation of condition (\ref{Cond Alpha seq}) is verified with $\gamma=1-2\alpha,$ and moreover, with   $\alpha\in (0,1/4),$ one has
  $$\sum_{k=0}^\infty \frac{\alpha_k\sqrt{\ln\beta_k}}{\beta_kA_k}=\sum_{k=0}^\infty O\left(\frac{a^{k^\alpha}k^{\alpha-1}\sqrt{\ln k^{1/2-\alpha}}}{k^{1/2-\alpha}a^{k^\alpha}}\right)=\sum_{k=0}^\infty O\left(\frac{\sqrt{\ln k}}{k^{3/2-2\alpha}}\right)<+\infty.$$
 	\vskip 0.1cm
 	
\subsection{Algorithms with the storage the past samples but updating subgradients}
In Algorithms 1 and 2, at the $k$-th iteration, we only need to compute a sugradient of the function $h(\cdot,s^k)$ at $x^k$. In this subsection, we propose the following two algorithms, in which all subgradients of the functions $h(\cdot,s^ i)$, $i=0,1,...,k$ are computed at the current iteration $x^k$. Let $\{\alpha_k\}$ be a sequence of positive reals with $\sum_{k=0}^\infty\alpha_k=+\infty.$
\vskip 0.1cm
\noindent\rule{12.7cm}{1.5pt}\\
\noindent \textbf{Algorithm 3:  Stochastic DC Algorithm 3 (SDCA3) }\\
\noindent\rule{12.7cm}{1pt}\\
 \texttt{Initialization:} Initial data: $x^0\in\dom r_1,$ draw $s^0\overset{\textrm{iid}}{\simeq} \mathbb{P},$  set $k=0.$ 
 \vskip 0.1cm
\texttt{Repeat:}
\begin{itemize}
 \item[1.] Compute $z^k_{i}\in \partial h(x^k,s^ i),\;\; i=0,1,...,k,$ and $w^k\in\partial r_2(x^k).$
 
 \item [2.] Set $y^k=\frac{1}{A_k}\sum_{i=0}^k\alpha_iz^k_i.$

\item [3.]  Compute a solution $x^{k+1}$ of the convex program
\begin{equation}\label{CP5}
\min \{G(x)+r_1(x)-\langle y^{k}+w^k,x\rangle :\ x\in \Bbb{R}^{n}\}.
\end{equation}
\item[4.] Set $k:=k+1$ and draw $s^{k}\overset{\textrm{iid}}{\simeq} \mathbb{P}.$
\end{itemize}

\texttt{Until} Stopping criterion.

\noindent\rule{12.7cm}{1.5pt}
\vskip 0.1cm

\bigskip

As SDCA3 implicitly supposes the form of $G$ is known explicitly, we further propose another scheme called SDCA4 to handle cases where $G$ is unknown by nature.

\noindent\rule{12.7cm}{1.5pt}\\
\noindent \textbf{Algorithm 4:  Stochastic DC Algorithm 4 (SDCA4) }\\
\noindent\rule{12.7cm}{1pt}

Similar to Algorithm 3, where $G$ in step 3 of Algorithm 3 is replaced by $G_k$,
\begin{align*}
G_k(x) = \frac{1}{A_k}\sum_{i=0}^k \alpha_ig(x,s^ i).
\end{align*}
\noindent\rule{12.7cm}{1.5pt}
\vskip 0.1cm
\noindent\textit{ Remark 2.} The principal difference between these algorithms  and Algorithms 1 and 2  is that at the $k$-th iteration, the subgradients of $h(\cdot,s^ i)$ ($i=1,...,k-1$) with respect to the past sample realizations are all updated. 
For Algorithms 3 and 4, we obtain a stronger convergence result that  almost surely the sequence $\{f(x^k)\}$ converges and all limit points of $\{x^k\}$ are critical points. Moreover, with the same added assumptions as for Algorithms 1 and 2, the convergence rate of $d_k$ is of order $O(\ln k/\sqrt{k}).$
\begin{theorem}\label{Convergence 2} Let a sequence of positive reals $\{\alpha_k\}$ be defined the same as Theorem \ref{General Convergence}. Suppose that (A1) holds for Algorithm 3 and (A2) holds for Algorithm 4. Let $\{x^k\}$ be a sequence generated by either Algorithm 3 or Algorithm 4. Suppose that the optimal value $\alpha$ of problem (\ref{SDCP}) is finite and with probability 1, $\limsup_{k\to\infty}\|x^k\|<\infty.$ Assuming $\rho:=\inf_{s\in\Omega} \rho(h(\cdot,s))+\rho(r_2)>0,$ then one has
\begin{itemize}
\item[(i)] Almost surely the sequence of function values $\{f(x^k)\}$ converges; $\sum_{k=0}^\infty\|x^{k+1}-x^k\|^2<+\infty,$ and all limit points of $\{x^k\}$ are critical points of (\ref{SDCP}).  
\item[(ii)] With the added assumptions as in Theorem \ref{General Convergence} $(ii)$, one has
$$d^2_k=O\left(\frac{1}{k}\sum_{i=1}^k\frac{(1+\sqrt{\ln \beta_{i-1}})^2}{\beta_{i-1}^2}\right).$$ In particular, for the sequence $\alpha_k=1,$ for all $k\in\N,$ then
$d_k=O(\ln k/\sqrt{k})$ as $k\to\infty.$  
\end{itemize} 
\end{theorem}
\vskip 0.1cm
\textit{Proof.} \emph{We only prove Theorem \ref{Convergence 2} for Algorithm 3. The proof for Algorithm 4 is analogous where we replace $G$ by $G_k$ and employ additional assumptions imposed on the function $g.$}

\textit{(i)} Consider the sequence $\{x^k\}$ generated by Algorithm 3. Firstly, we denote by 
$$\mathcal{F}_k=\sigma(s^0,...,s^{k-1}, x^0,...,x^{k},z^0,...,z^{k-1},u^0,...,u^{k-1}),\quad k\in\mathbb{N},$$
the increasing $\sigma-$ field generated by random variables $s^0,...,s^{k-1},x^0,...,x^{k},z^0,...,z^{k-1},$ $u^0,...,u^{k-1}.$ For $k=1,2,...,$ define the function $V_k:\R^n\to\R\cup\{+\infty\}$, $k\in\mathbb{N}$,
$$V_k(x)=G(x)+r_1(x)-\left\langle\frac{1}{A_k}\sum_{i=0}^{k} \alpha_iz^k_i+w^k,x-x^{k}\right\rangle-\frac{1}{A_k}\sum_{i=0}^{k}\alpha_ih(x^k,s^ i)-r_2(x^k).$$
As $x^{k+1}$ is a solution of the problem (\ref{CP5}), one has
\begin{equation}\label{estim 5}
V_{k}(x^{k+1})\le V_k(x^k)=G(x^k)+r_1(x^k)-\frac{1}{A_k}\sum_{i=0}^k\alpha_ih(x^k,s^ i)-r_2(x^k).
\end{equation}
By the strong convexity of the functions $h(\cdot,s)+r_2,$ $s\in\Omega$ with modulus at least $\rho,$
\begin{equation}\label{estim 5 bis}
\begin{array}{lll}
&V_{k-1}(x^k)-V_k(x^k)=\frac{1}{A_k}\sum_{i=0}^k\alpha_ih(x^k,s^ i)\\
&-\frac{1}{A_{k-1}}\left[\sum_{i=0}^{k-1}\alpha_i\left(\langle z^{k-1}_i,x^k-x^{k-1}\rangle + h(x^{k-1},s^ i)\right)\right] \\
&+[r_2(x^k)-\langle w^{k-1},x^k-x^{k-1}\rangle - r_2(x^{k-1})]\\
&\ge \frac{\alpha_k}{A_k}\left[h(x^k,s^k)-\frac{1}{A_{k-1}}\sum_{i=0}^{k-1}\alpha_ih(x^k,s^ i)\right]+\frac{\rho}{2}\|x^k-x^{k-1}\|^2.
 \end{array}
\end{equation}
Thus 
\begin{equation}\label{Decreasing-V}
\begin{array}{ll}
V_{k}(x^{k+1})&\le V_{k-1}(x^{k})-(V_{k-1}(x^k)-V_k(x^k))\\
&\le V_{k-1}(x^{k})-\frac{\alpha_k}{A_k}\left[h(x^k,s^k)-\frac{1}{A_{k-1}}\sum_{i=0}^{k-1}\alpha_ih(x^k,s^ i)\right]-\frac{\rho}{2}\|x^k-x^{k-1}\|^2. 
\end{array}
\end{equation}
By taking expectation with respect to $\mathcal{F}_k$ both sides, we derive 
\begin{equation}\label{Estim-Expec-bis}
\begin{array}{ll}
\mathbb{E}_{\mathcal{F}_k}V_{k}(x^{k+1})&\le V_{k-1}(x^{k})-\frac{\alpha_k}{A_k}\left[H(x^k)-\frac{1}{A_{k-1}}\sum_{i=0}^{k-1} \alpha_ih(x^{k},s^ i)\right]-\frac{\rho}{2}\|x^k-x^{k-1}\|^2\\
&\le V_{k-1}(x^{k})-\xi_k-\frac{\rho}{2}\|x^k-x^{k-1}\|^2,
\end{array}
\end{equation}
where $\xi_k=\frac{\alpha_k}{A_k}\left[H(x^k)-\frac{1}{A_{k-1}}\sum_{i=0}^{k-1}\alpha_ih(x^k,s^ i)\right].$
By Lemma \ref{ULLN}, for some $c>0,$
\begin{equation}\label{Ge expec xi-bis}
\mathbb{E}|\xi_k|\le\frac{c\alpha_k(1+\sqrt{\ln \beta_{k-1}})}{\beta_{k-1}A_k},\;\;\mbox{for all}\; k= 1,2,...
\end{equation}
 Hence, $\sum_{k=1}^\infty\mathbb{E}|\xi_k|<+\infty.$ Let $R>0$ such that  the sequence $\{x^k\}$ is contained in a compact set $C\subseteq B_R=\{x\in\R^n;\:\|x\|<R\}.$ In view of Claim 1 in the proof of the preceding theorem, with probability 1, we has
 \begin{equation}\label{Conver Prob 1}
 \lim_{k\to\infty}\frac{1}{A_k}\sum_{i=0}^k\alpha_ih(x,s^i)=H(x),\;\;\forall x\in\overline{B}_R.
 \end{equation}
 From relations (\ref{estim 5}), (\ref{estim 5 bis}), one has the following estimates
 \begin{equation}\label{estim @@@}
 \begin{array}{lll}
 &V_k(x^{k+1})\le V_k(x^k)= f(x^k)+ \left[ H(x^k)-\frac{1}{A_k}\sum_{i=0}^{k}\alpha_ih(x^k,s^ i)\right]\\
&\le V_{k-1}(x^k)- \frac{\alpha_k}{A_k}\left[h(x^k,s^k)-\frac{1}{A_{k-1}}\sum_{i=0}^{k-1}\alpha_ih(x^k,s^ i)\right]-\frac{\rho}{2}\|x^k-x^{k-1}\|^2.\\
 \end{array}
 \end{equation}
 As $\alpha_k/A_k\to 0$ (by the first condition of (\ref{Cond Alpha seq})) and $h(\cdot,s)$ is uniformly bounded on $\overline{B}_R,$ in view of (\ref{Conver Prob 1}), this relation implies that $\{V_k(x^{k+1})\}$ is bounded below.  Thus, thanks to the supermartingale convergence, we  arrive at the conclusion that almost surely the sequence $\{V_k(x^{k+1})\}$ converges and $\sum_{k=0}^\infty\|x^{k+1}-x^k\|^2<+\infty.$ By relation (\ref{estim @@@}), the convergence of $\{V_k(x^{k+1})\}$ yields immediately the convergence of $\{V_k(x^k)\}$ and $\{f(x^k)\}$.
 We denote 
 \begin{equation}\label{ Event almost sure}
 S:=\{\{s^k\}_{k\in\N}:\;\;\lim_{k\to\infty}\frac{1}{A_k}\sum_{i=0}^k\alpha_ih(x,s^i)=H(x),\;\|x^k-x^{k-1}\|\to0\;\;\forall x\in\overline{B}_R \}.
 \end{equation}
 Then by (\ref{Almost surely}), $\mathbb{P}(S)=1.$ For given $\{s^k\}_{k\in\N}\in S,$ and suppose that $\{x^k\}_{k\in\N}$ is a sequence generated by Algorithm 3 with respect to the sequence of samples $\{s^k\}.$ For any limit point  $x^*$ of $\{x^k\},$ pick a subsequence $\{x^{l_k}\}$ converging to $x^*.$ As $\|x^k-x^{k-1}\|\to 0,$ then $x^{l_k-1}\to x^*.$ By passing to a subsequence, we can assume that $\lim_{k\to\infty}y^{l_k-1}=y^*$ and $\lim_{k\to\infty}w^{l_k-1}=w^*.$ As $y^{l_k-1}+w^{l_k-1}\in \partial G(x^{l_k})+\partial r_1(x^{l_k}),$ and $w^{l_k-1}\in\partial r_2(x^{l_k-1}),$  one has $y^*+w^*\in\partial G(x^*)+\partial r_1(x^*)$ and $w^*\in\partial r_2(x^*).$ On the other hand, as
 $y^{l_k-1}\in \frac{1}{A_{l_k-1}}\sum_{i=0}^{l_k-1}\alpha_i\partial h(x^{l_k-1},s^i),$
 \begin{equation}\label{Subdiff inq}
 \langle y^{l_k-1}, x-x^{l_k-1}\rangle \le \frac{1}{A_{l_k-1}}\sum_{i=0}^{l_k-1}\alpha_ih(x,s^i)-\frac{1}{A_{l_k-1}}\sum_{i=0}^{l_k-1}\alpha_ih(x^{l_k-1},s^i)\;\;\forall x\in\overline{B}_R. 
 \end{equation}
 Note that $\frac{1}{A_{l_k-1}}\sum_{i=0}^{l_k-1}\alpha_ih(x,s^i)\to H(x)$ for all $x\in\overline{B}_R,$ and $\frac{1}{A_{l_k-1}}\sum_{i=0}^{l_k-1}\alpha_ih(x^{l_k-1},s^i)\to H(x^*),$ where the later relation follows from
 $$\left|\frac{1}{A_{l_k-1}}\sum_{i=0}^{l_k-1}\alpha_ih(x^{l_k-1},s^i)-H(x^*)\right|\le \left|\frac{1}{A_{l_k-1}}\sum_{i=0}^{l_k-1}\alpha_ih(x^*,s^i)-H(x^*)\right|+ L\|x^{l_k-1}-x^*\|,$$
 here $L$ is some common Lipschitz constant of $h(\cdot,s)$ on $\overline{B}_R.$ Hence, by letting $k\to\infty$ in (\ref{Subdiff inq}), one obtains
 $\langle y^*,x-x^*\rangle\le H(x)-H(x^*),\;\;\forall x\in \overline{B}_R.$
 Since $x^*\in C\subset B_R,$ the inequality yields $y^*\in\partial H(x^*),$ showing that
 $x^*$ is a critical point of problem (\ref{SDCP}).  
 \vskip 0.1cm
\textit{(ii)} Assume now $\nabla r_2$ and the derivatives $\nabla h(\cdot,s),$ $s\in \Omega,$ are equi-Lipschitz with modulus $L/2$ on a compact set containing $\{x^k\}.$ For $k\in\N_*,$
\begin{equation}\label{Deriv-Lip-bis}
\begin{array}{ll}
y^{k-1}+w^{k-1}&= \frac{1}{A_{k-1}}\sum_{i=0}^{k-1}\alpha_i\nabla h(x^{k-1},s^ i)+\nabla r_2(x^{k-1})\\
&\in \frac{1}{A_{k-1}}\sum_{i=0}^{k-1}\alpha_i\nabla h(x^k,s^ i)+\nabla r_2(x^k)+L\|x^k-x^{k-1}\|\mathbb{B}\\
&\subseteq\nabla H(x^k)+\nabla r_2(x^k) +\left(\eta_k +L\|x^k-x^{k-1}\|\right)\mathbb{B},
\end{array}
\end{equation}
where, $\eta_k=\left\|\frac{1}{A_{k-1}}\sum_{i=0}^{k-1}\alpha_i\nabla h(x^k,s^ i)-\nabla H(x^k)\right\|.$
Furthermore, by noticing that $y^{k-1}+w^{k-1}\in\partial G(x^k)+\partial r_1(x^k),$ one has
\begin{equation}\label{dist-Critical-bis}
d(0,\partial G(x^k)+\partial r_1(x^k)-\nabla H(x^k)-\nabla r_2(x^k))\le \eta_k+L\|x^k-x^{k-1}\|. 
\end{equation}
 By Lemma \ref{ULLN}, for some $c>0,$
\begin{equation}\label{Inq-Dist-bis}
\tau_k:=\mathbb{E}d(0,\partial G(x^k)+\partial r_1(x^k)-\nabla H(x^k)-\nabla r_2(x^k))\le L \mathbb{E}\|x^k-x^{k-1}\|+ \frac{c(1+\sqrt{\ln \beta_{k-1}})}{\beta_{k-1}}.
\end{equation}
As $\mathbb{E}\|x^k-x^{k-1}\|^2\ge \tau_k^2/(2L^2)-\frac{c^2(1+\sqrt{\ln \beta_{k-1}})^2}{L^2\beta_{k-1}^2},$ and by (\ref{Ge expec xi-bis}) and (\ref{Decreasing-V}), 
$$\begin{array}{ll}\mathbb{E}V_{k}(x^{k+1})- \mathbb{E}V_{k-1}(x^k)&\le -\frac{\rho}{2}\mathbb{E}\|x^k-x^{k-1}\|^2+\mathbb{E}|\xi_k|\\
&\le -\frac{\rho\tau_k^2}{4L^2}+\frac{\rho c^2(1+\sqrt{\ln \beta_{k-1}})^2}{2L^2\beta_{k-1}^2} +\frac{c\alpha_k(1+\sqrt{\ln \beta_{k-1}})}{\beta_{k-1}A_k}.
\end{array}$$ 
 Therefore, for $k\in\N_*,$
$$
\frac{\rho}{4L^2}\sum_{i=1}^k \tau_i^2\le \mathbb{E}V_0(x^1)-\mathbb{E}V_{k}(x^{k+1})+\sum_{i=1}^k\frac{\rho c^2(1+\sqrt{\ln \beta_{i-1}})^2}{2L^2\beta_{i-1}^2}+\sum_{i=1}^k\frac{c\alpha_i(1+\sqrt{\ln \beta_{i-1}})}{\beta_{i-1}A_i},
$$
implying, for some $V^*>0,$ 
$$\frac{\rho}{4L^2}d^2_k\le \frac{\mathbb{E}V_0(x^1)-V^*+\sum_{i=1}^k\frac{\rho c^2(1+\sqrt{\ln \beta_{i-1}})^2}{2L^2\beta_{i-1}^2}+\sum_{i=1}^k
\frac{c\alpha_i(1+\sqrt{\ln \beta_{i-1}})}{\beta_{i-1}A_i}}{k}.$$ 
By the assumption
$\sum_{i=1}^\infty\frac{\alpha_i(1+\sqrt{\ln \beta_{i-1}})}{\beta_{i-1}A_i}<+\infty,$
we conclude $$d^2_k=O\left(\frac{1}{k}\sum_{i=1}^k\frac{(1+\sqrt{\ln \beta_{i-1}})^2}{\beta_{i-1}^2}\right).$$ Finally, when $\alpha_k=1,$ $k\in\N,$
$\sum_{i=1}^k\frac{(1+\sqrt{\ln \beta_{i-1}})^2}{\beta_{i-1}^2}=\sum_{i=1}^k\frac{(1+\sqrt{\ln i})^2}{i^2}=O(\ln^2k).$
 \hfill {$\Box$}

\section{Numerical Experiments}

\bigskip

Principal component analysis is one of the most successful tools for dimensionality reduction. In this section, instead of estimating the direction with the highest variability on given datasets, our aim is to generalize this performance on unseen data, that is, we consider the Expected problem of Principal Component Analysis (denoted by (E-PCA)) \cite{Montanari},
\begin{align*}
&\text{minimize} \quad \varPhi(x)= -\dfrac{1}{2} \mathbb{E}_{s \sim \mathbb{P}}(\langle x,s \rangle^2) d \mathbb{P}, \quad \text{(E-PCA)}\\
&\text{subject to} \quad \Vert x \Vert \leq 1,
\end{align*}
where $s$ is a normalized random vector, i.e., $\Vert s \Vert =1$ whose distribution is unknown.

We consider the situation where the training data is given on stream; hence we do not have the whole training dataset at the beginning. The coming data is then fed to SDCA schemes until the end of the stream. We perform experiments on standard datasets of LIBSVM \footnote{The datasets can be downloaded at \url{https://www.csie.ntu.edu.tw/~cjlin/libsvm/}}, namely \texttt{protein, YearPredictionMSD, SensIT Vehicle, shuttle}. Samples of each dataset are normalized as $\Vert s_i \Vert = 1.$ The training data is provided on stream, and each algorithm terminates when a dataset is used up. For each run, we randomly shuffle training datasets before delivering them to algorithms. Besides, the starting points are randomly initialized in $S$, where $S = \{x: \Vert x \Vert \leq 1\}$. We run each algorithm $20$ times and report the (average of $20$ runs) objective on validation sets along the execution process, which represents the expected objective (generalized performance). To enhance visualization, we plot the gap between the expected objective and the ``optimal value" found of the PCA problem constructed using validation data by running deterministic DCA  with multiple initial points ($10$ points randomly chosen in $S$).

All experiments are performed on a PC Intel(R) Core(TM) i7-8700 CPU@3.20GHz of 16 GB RAM.

To investigate different aspects of our algorithms, we propose three experiments. In the first experiment, we formulate the (E-PCA) problem as a stochastic DC program with
$ r_1(x) = \chi_S(x), r_2(x) = 0, g(x,s) = \dfrac{\lambda}{2} \Vert x \Vert^2, h(x,s) = \dfrac{\lambda}{2} \Vert x \Vert^2 + \dfrac{1}{2} \langle x,s \rangle^2,$
where $\lambda >0$ and $\chi_S$ is the indicator function of $S$. The positive hyperparameter $\lambda$ is used to guarantee our strong convexity assumption. In practice, we observed that small $\lambda$ yields better results, we therefore set $\lambda=10^{-6}$. With this expression, since $g(x,s)$ is independent of $s$, SDCA1 coincides with SDCA2 and SDCA2 coincides with SDCA4. It is clear that with this DC decomposition, assumption (A1) holds. Furthermore, it is noteworthy that, with this formulation, SDCA1 coincides with a version of SSUM \cite{Razaviyayn} with the  setting $g_2 = g,g_1 = -h, \mathcal{X} = S$ and $\hat{g}_1(x,y,\xi) = g_1(y,\xi) + \langle \nabla g_1(y,\xi),x-y \rangle$ (with the same notations used in \cite{Razaviyayn}). We compare our algorithms with Projected Stochastic Subgradient method (PSS) with the projected region being $S$. For PSS, we use two types of stepsize policies: constant and diminishing. The constant stepsize is searched in $\{0.001,0.005,0.01,0.015,0.02\}$. We found that $0.005$ is the most appropriate candidate. On the other hand, the sequence of diminishing stepsize is given by $\alpha_k = c/k$, where $c$ is another hyperparameter being searched in $\{4,5,\ldots,11,12\}$. A very good candidate found is $c = 8.$

\begin{figure}
     \centering
        \subfigure[\texttt{SensIT Vehicle}]{	
	\includegraphics[width=0.24\textwidth]{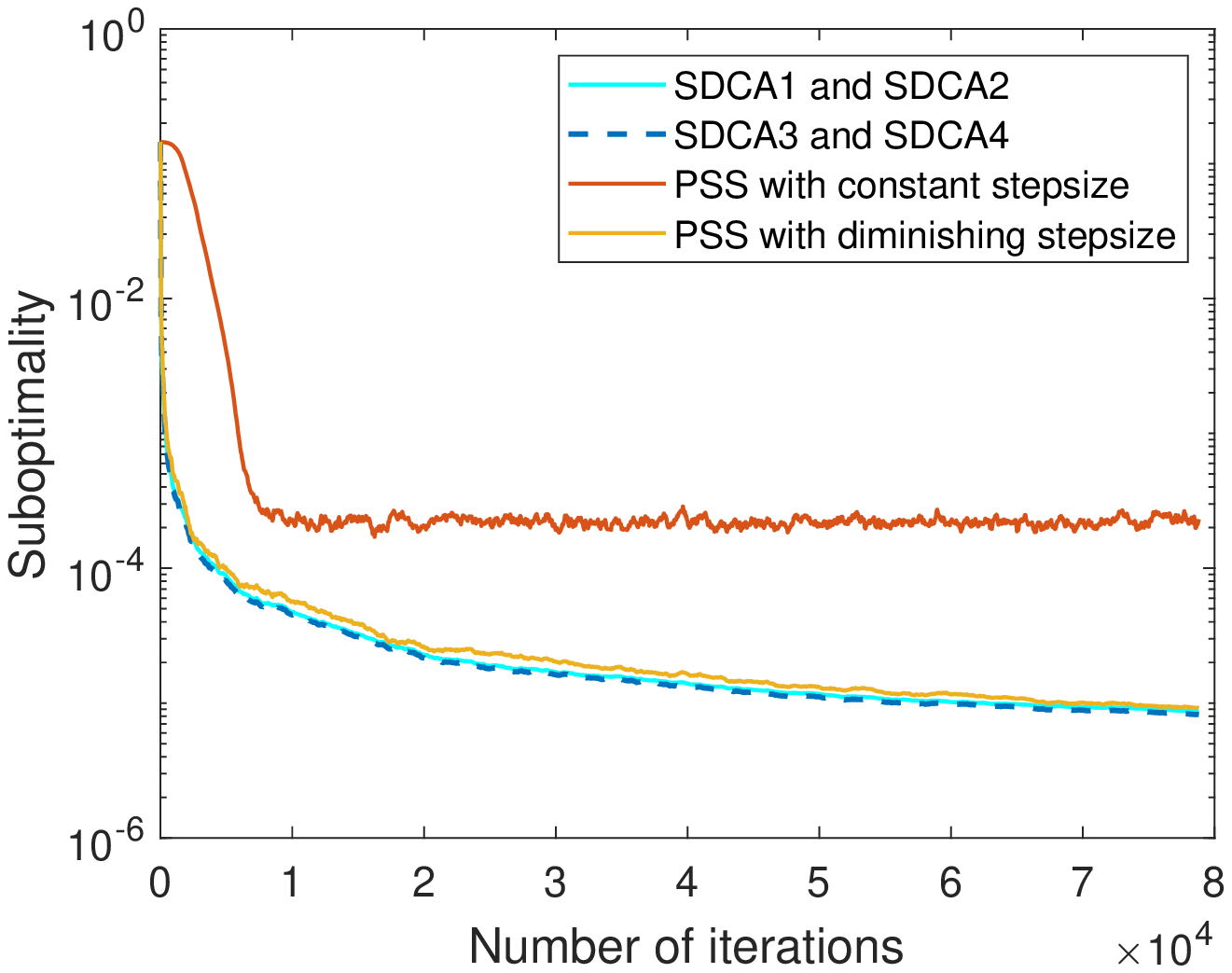} 
	}
	\hspace{-10pt}
	\subfigure[\texttt{shuttle}]{	
	\includegraphics[width=0.24\textwidth]{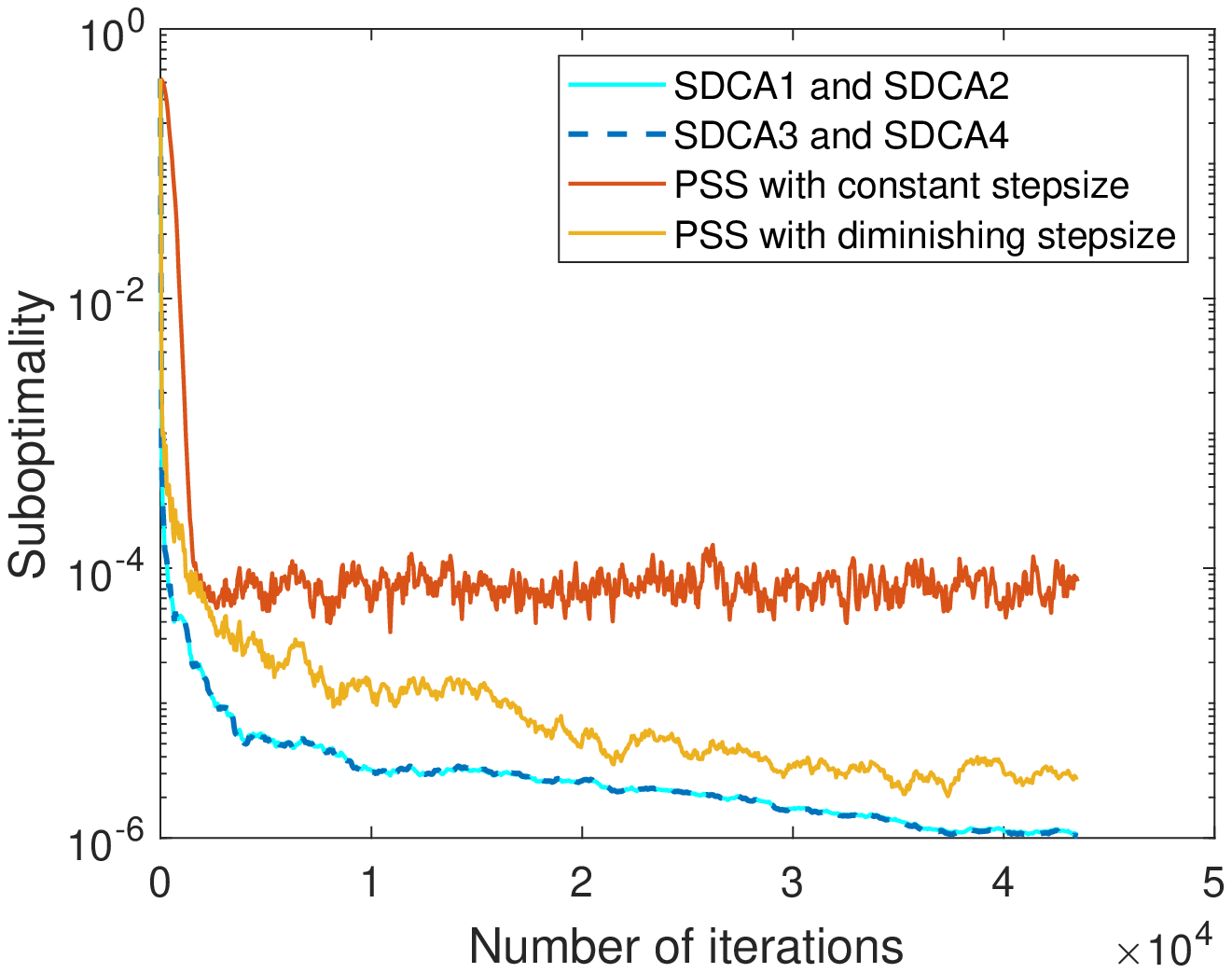} 
	}
	\hspace{-10pt}
	\subfigure[\texttt{protein}]{	
	\includegraphics[width=0.24\textwidth]{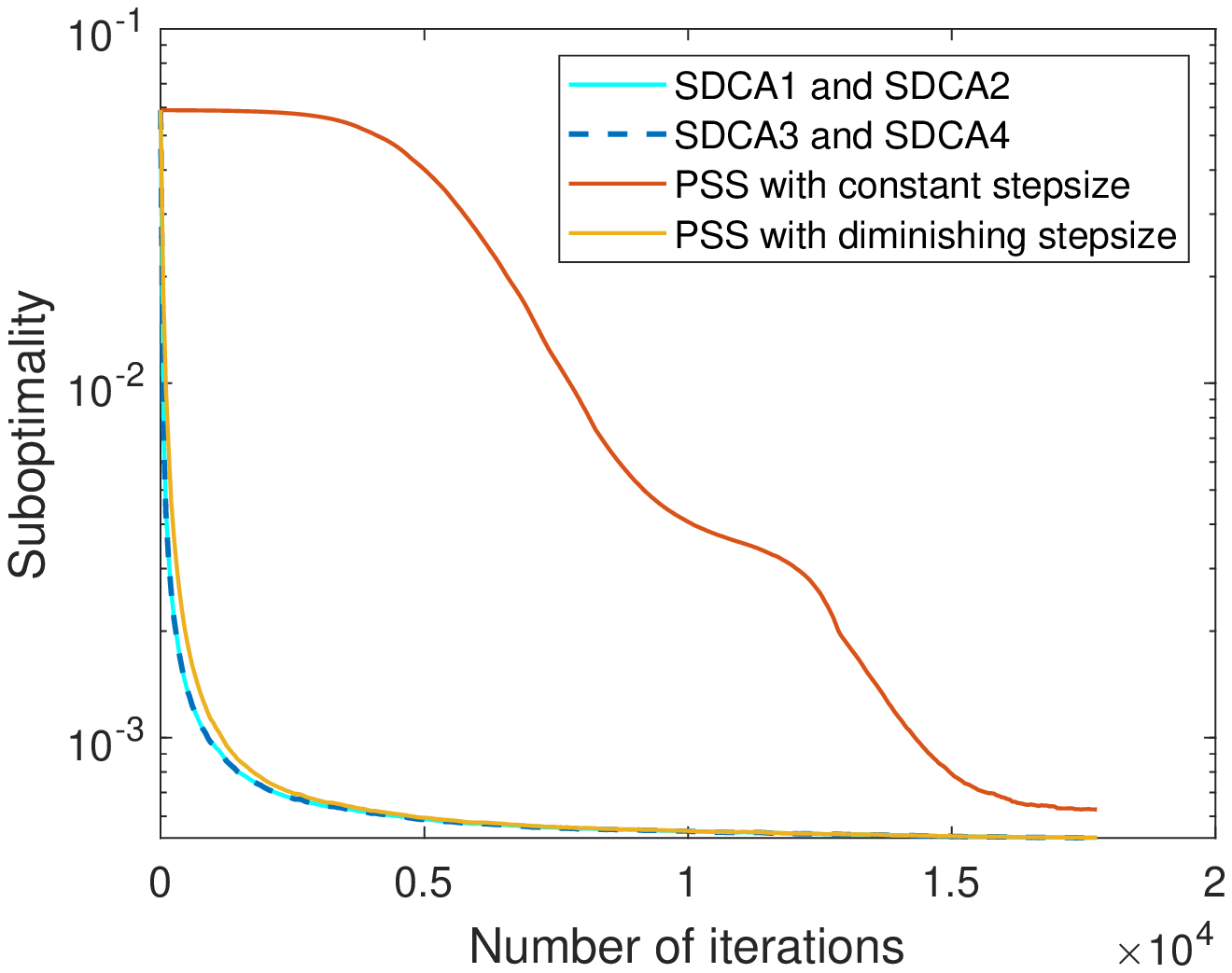} 
	}
	\hspace{-10pt}
	\subfigure[\texttt{YearPredictionMSD}]{	
	\includegraphics[width=0.24\textwidth]{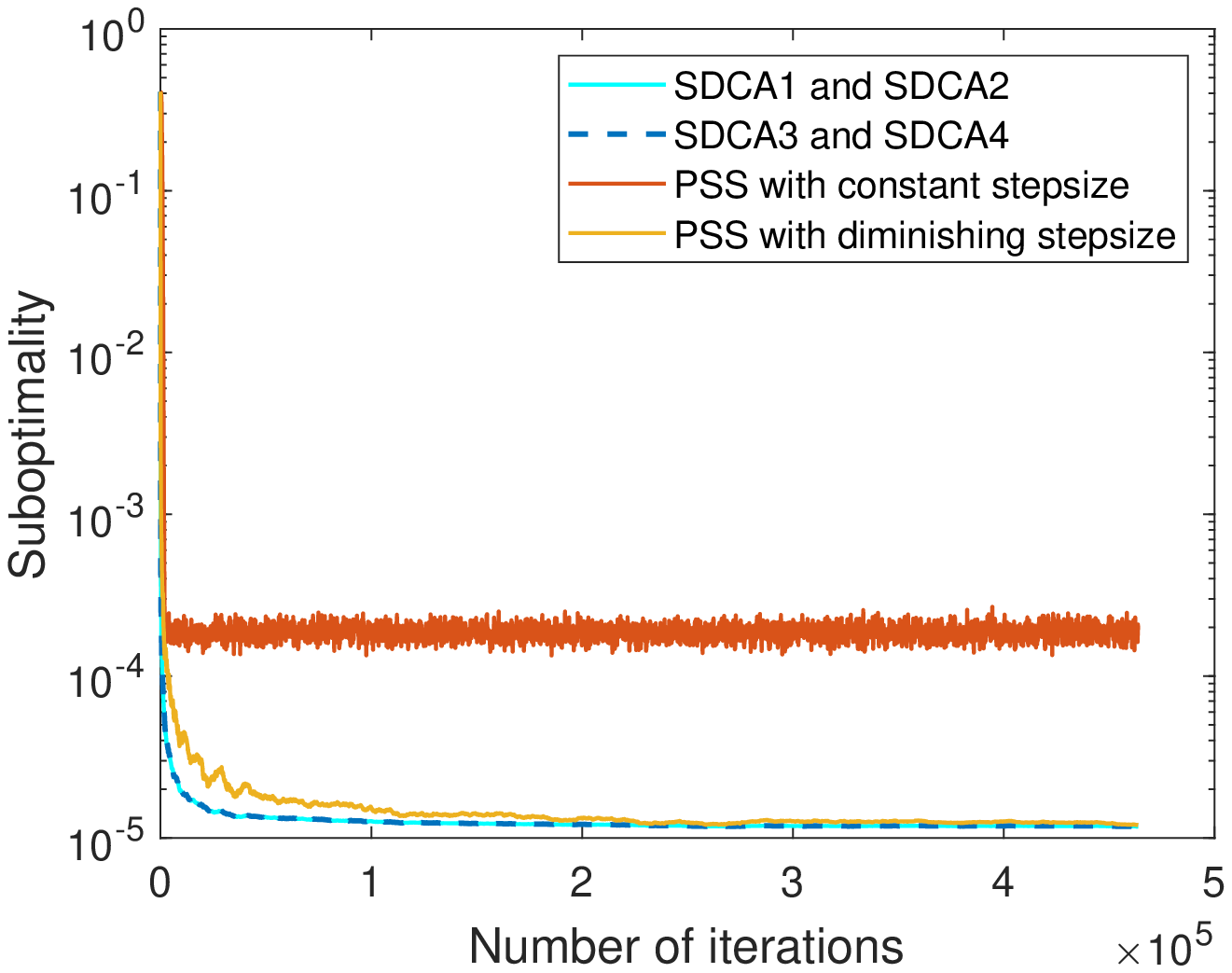} 
	}
	\hspace{-10pt}
	\subfigure[\texttt{SensIT Vehicle}]{	
	\includegraphics[width=0.24\textwidth]{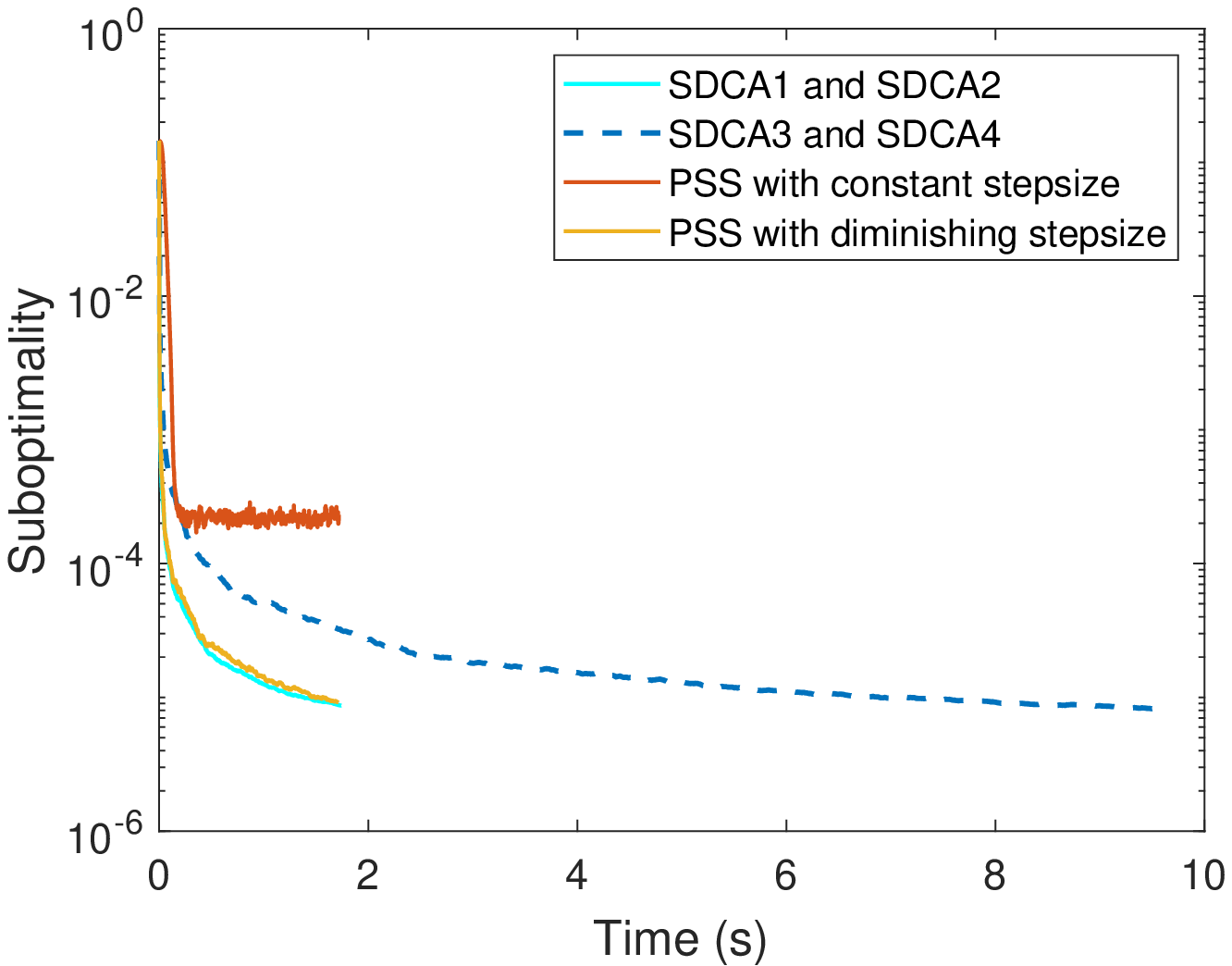} 
	}
	\hspace{-10pt}
	\subfigure[\texttt{shuttle}]{	
	\includegraphics[width=0.24\textwidth]{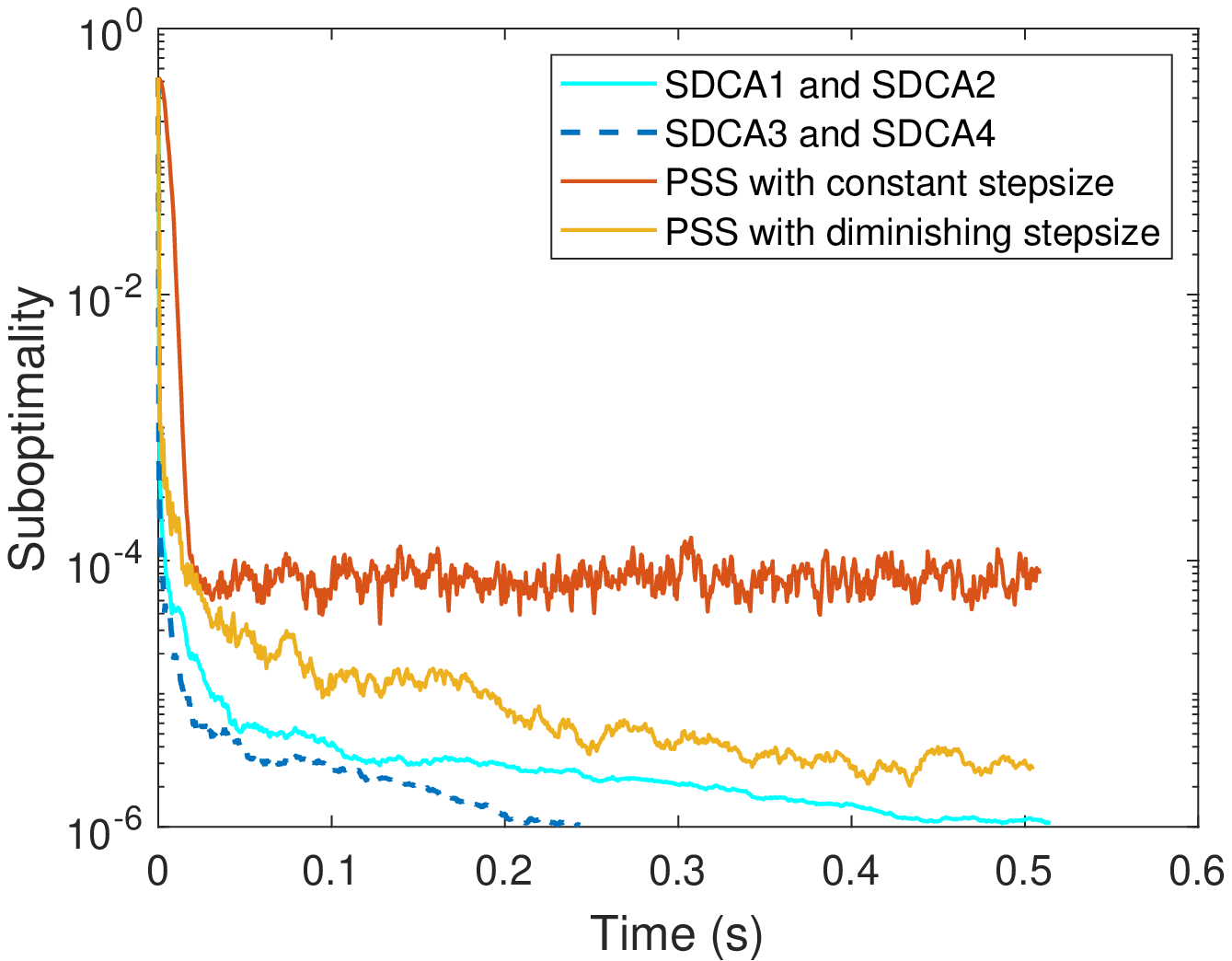} 
	}
	\hspace{-10pt}
	\subfigure[\texttt{protein}]{	
	\includegraphics[width=0.24\textwidth]{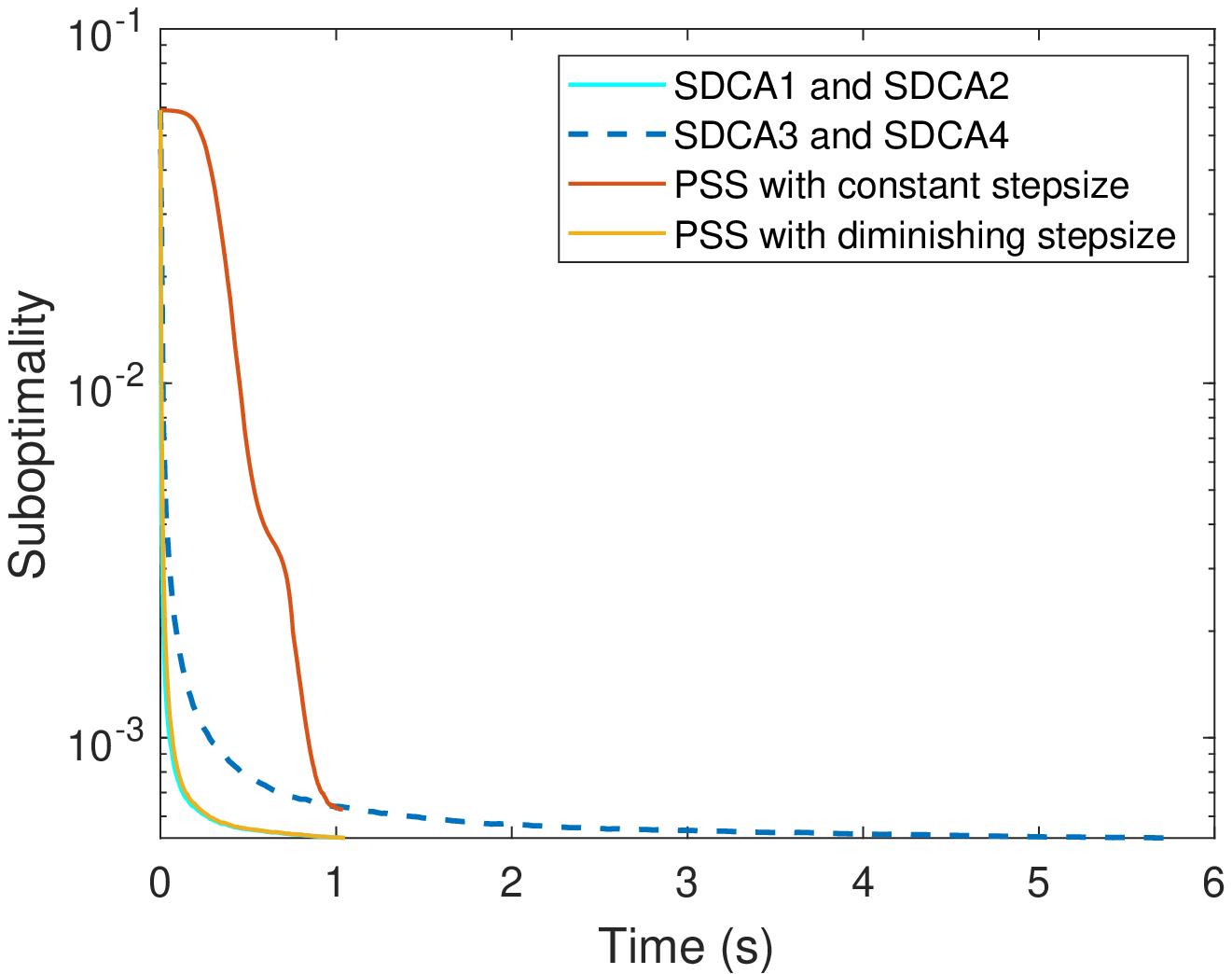} 
	}
	\hspace{-10pt}
	\subfigure[\texttt{YearPredictionMSD}]{	
	\includegraphics[width=0.24\textwidth]{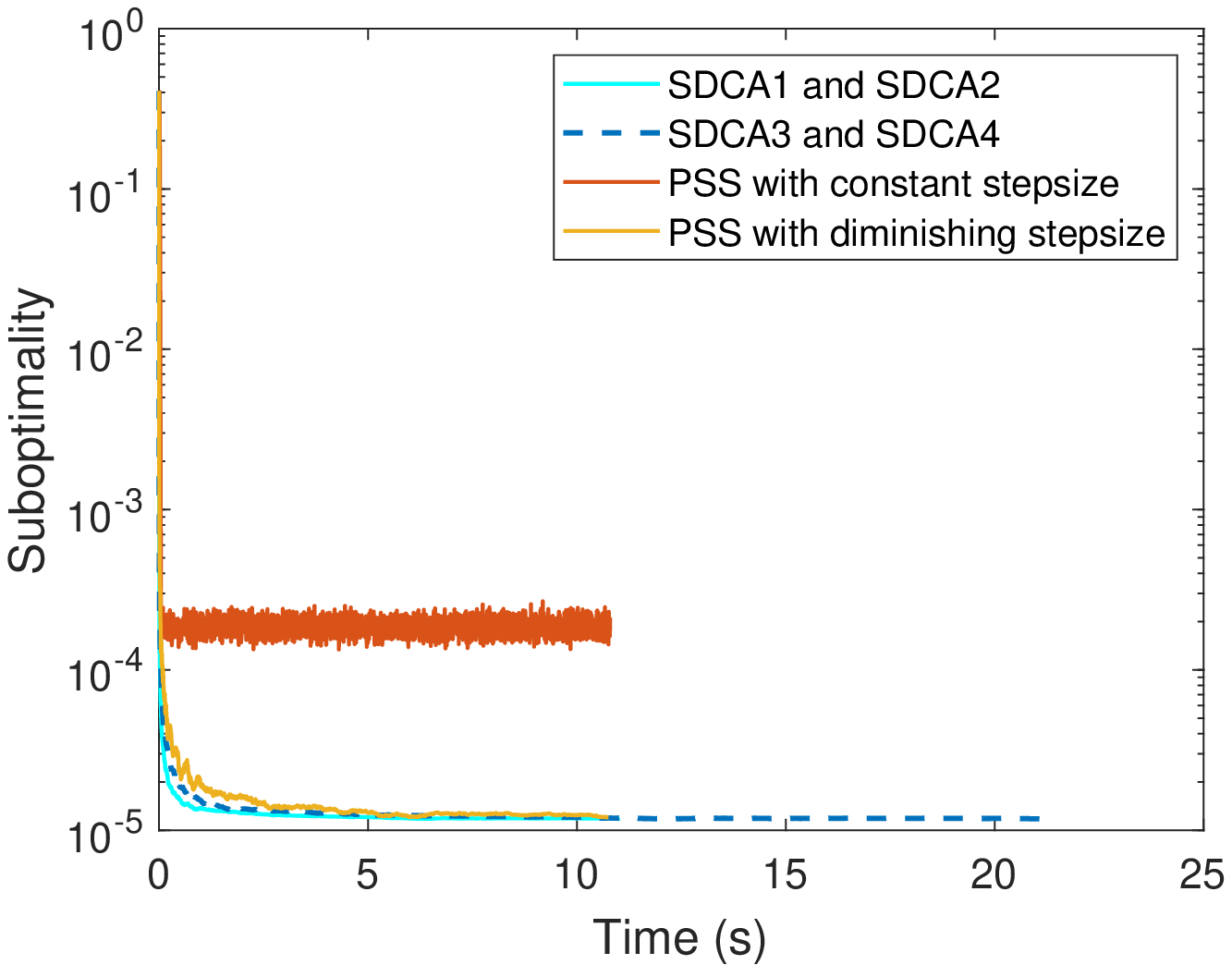} 
	}

        \caption{The performance of SDCA schemes compared with PSS methods}
        \label{fig1}
\end{figure}

We report the evolution of the objective on validation sets where the horizontal bar counts the number of iterations. It is noted that, in each iteration, each of our algorithms as well as PSS consume one new fresh sample (SDCA3 and SDCA4 use one new sample and reuse old samples); therefore, the cost of sample retrieving for one iteration is the same for all algorithms. Moreover, each iteration of these algorithms also requires solving one convex sub-problem. The first row of figure \ref{fig1} presents the performance of these algorithms in this regard. We observe that the performances of SDCA1,2 and SDCA3,4 are almost identical and better than the performance of two versions of PSS. Our algorithms achieve very good objective values where the suboptimality (considered at the end of the process) ranges from $1.06 \times 10^{-6}$ to $5.22 \times 10^{-4}$ (SDCA1,2) and from $1.03 \times 10^{-6}$ to $5.21 \times 10^{-4}$ (SDCA3,4). In all datasets, the performance of SDCA3,4 is a bit better than SDCA1,2 where the gap (the difference of two objective at the end of the process) varies from $8.93 \times 10^{-10}$ to $4.72 \times 10^{-7}$. For the PSS algorithms, while PSS with constant stepsize struggles to approach the optimal value, the diminishing stepsize policy obtains better performance as it almost reaches the performance of SDCA schemes. Nevertheless, in four datasets, PSS with diminishing stepsize is still a bit inferior to SDCA schemes where the gap (compared with SDCA1,2) varies from $2.46 \times 10^{-7}$ to $1.77 \times 10^{-6}$. 

On the other hand, we observe that the per-iteration stochastic-gradient-computing complexity (based on the number of stochastic gradients computed at each iteration) of PSS and SDCA1,2 is the same which is $\mathcal{O}(1)$, while the per-iteration stochastic-gradient-computing complexity of SDCA3,4 is $\mathcal{O}(k)$. Therefore, we further plot the evolution of the objective along the execution time horizon to study the performance of these algorithms regarding computational cost (figure \ref{fig1}, the second row). As the figure well illustrates, while the execution time of SDCA1,2 and two versions of PSS is similar, SDCA3,4 need more time to proceed through training sets (except for the dataset \texttt{shuttle}). The ratios of execution time of SDCA3,4 over SDCA1,2 are $5.48, 0.47, 5.55, 1.94$ on \texttt{sensIT Vehicle, shuttle, protein, YearPredictionMSD}, respectively. Since the objective gain of SDCA3,4 compared to SDCA1,2 is negligible while the running time is considerably longer, SDCA1,2 would be a better choice in this experiment.

In the second experiment, we want to study a ``really" stochastic program where $G$ is unknown by nature, and therefore SDCA1 and SDCA3 fail to work. In such a situation, our aim is to guarantee SDCA2 and SDCA4 still perform well. We first observe that for each $s, \Vert s \Vert =1,$ $\varphi(\cdot,s)$ is $1$-smooth on $S$, where $\varphi(x,s) = -\frac{1}{2} \langle x,s \rangle^2$  Consequently, $\frac{L}{2} \Vert x \Vert^2 \pm \varphi(x,s)$ are convex on $S$ if $L \geq 1$. Therefore, we have another DC reformulation for the (E-PCA) problem as follows
\begin{align*}
r_1(x) = \chi_S(x), r_2(x) = 0,g(x,s) = \dfrac{L}{2} \Vert x \Vert^2 - \dfrac{1}{2} \langle x,s\rangle^2, h(x,s) = \dfrac{L}{2} \Vert x \Vert^2 + \dfrac{1}{2} \langle x,s \rangle^2.
\end{align*}
With this setting, SDCA1 and SDCA3 are no longer applicable. Moreover, it is easy to verify that assumption $(A2)$ holds in this case. We choose $L=1.1$ which is a neutral parameter. It is worth mentioning that SDCA2 in this case coincides with a version of SSUM \cite{Razaviyayn} with the  setting $g_2 = g,g_1 = -h, \mathcal{X} = S$ and $\hat{g}_1(x,y,\xi) = g_1(y,\xi) + \langle \nabla g_1(y,\xi),x-y \rangle$. At iteration $k$, SDCA2 and SDCA4 require minimizing the following convex function 
\begin{align*}
    \tilde{f}(x)=\underbrace{\dfrac{L}{2} \Vert x \Vert^2 -\langle y^k, x \rangle}_{\tilde{g}(x)} - \underbrace{\dfrac{1}{2(k+1)} \sum_{i=0}^{k}{\langle x,s^i \rangle^2}}_{\tilde{h}(x)}.
\end{align*}

\begin{figure}
     \centering
        \subfigure[\texttt{SensIT Vehicle}]{	
	\includegraphics[width=0.24\textwidth]{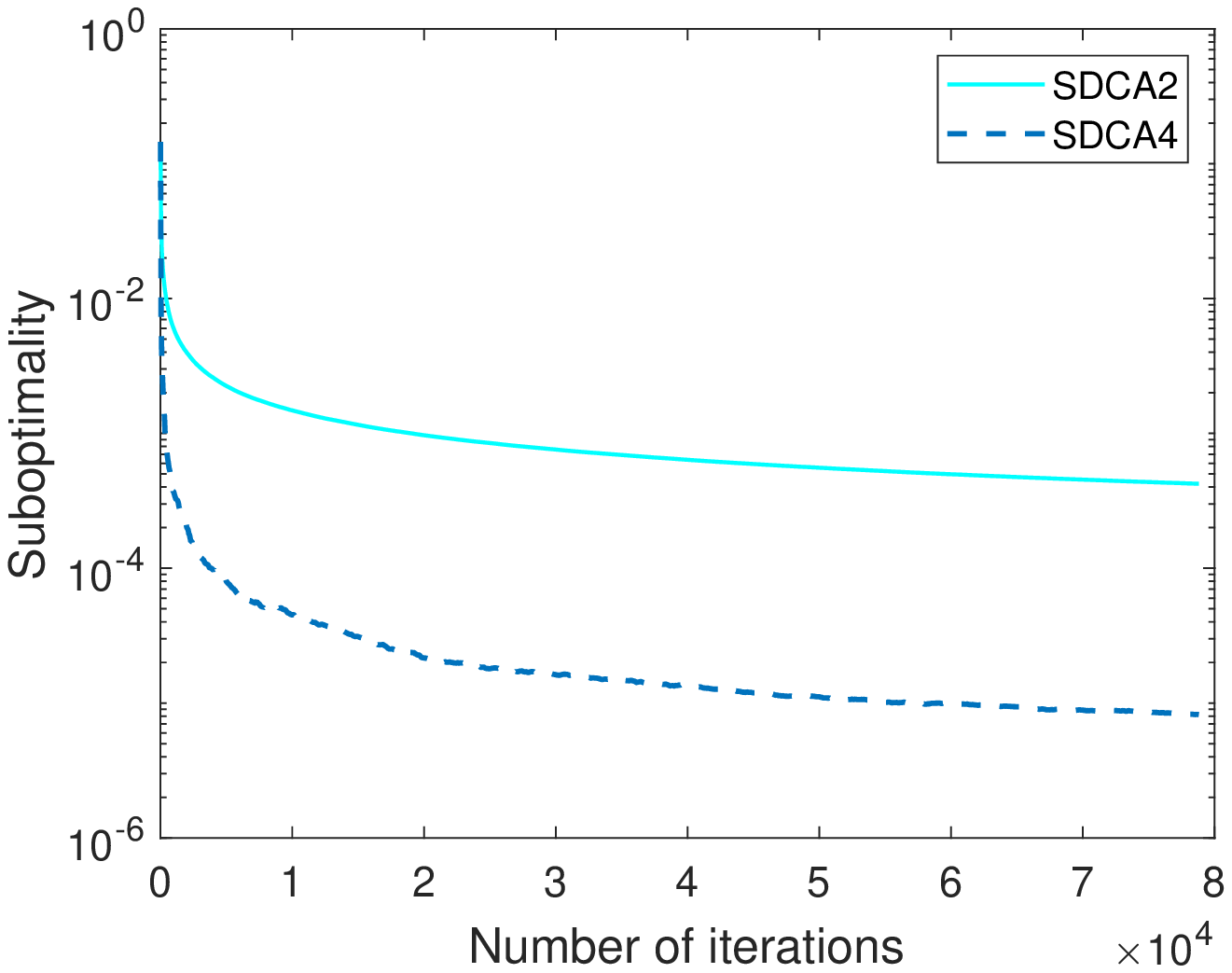} 
	}
	\hspace{-10pt}
	\subfigure[\texttt{shuttle}]{	
	\includegraphics[width=0.24\textwidth]{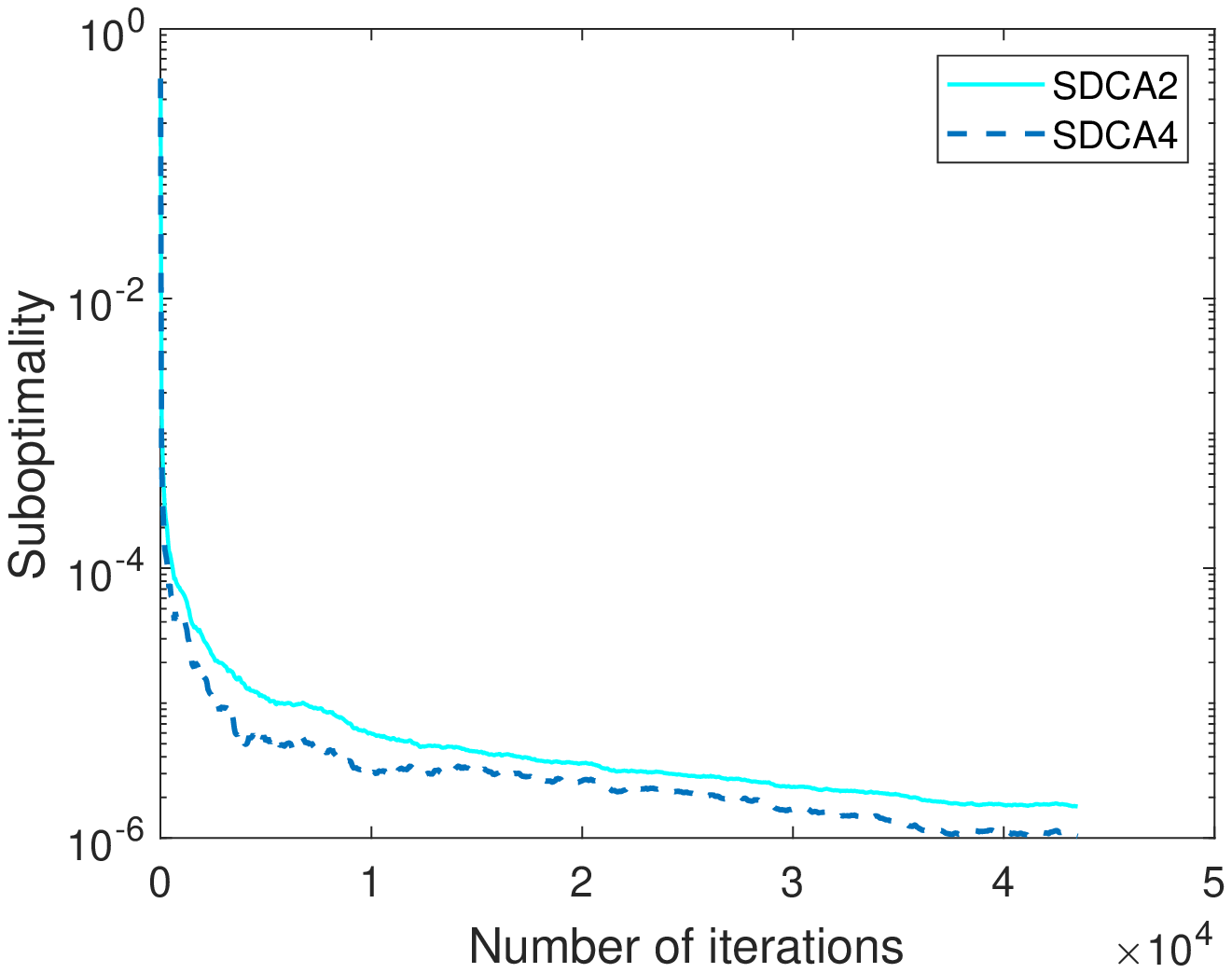} 
	}
	\hspace{-10pt}
	\subfigure[\texttt{protein}]{	
	\includegraphics[width=0.24\textwidth]{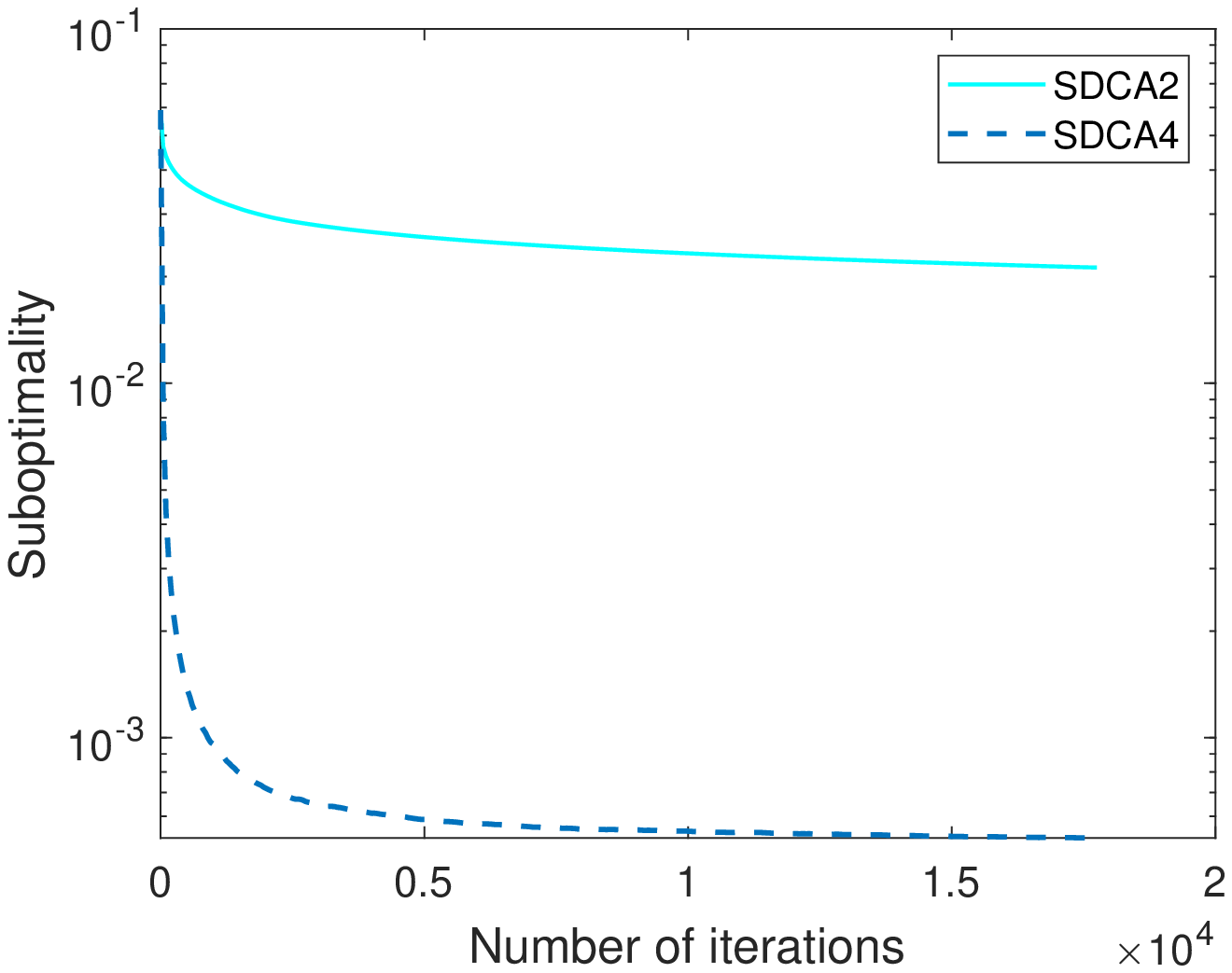} 
	}
	\hspace{-10pt}
	\subfigure[\texttt{YearPredictionMSD}]{	
	\includegraphics[width=0.24\textwidth]{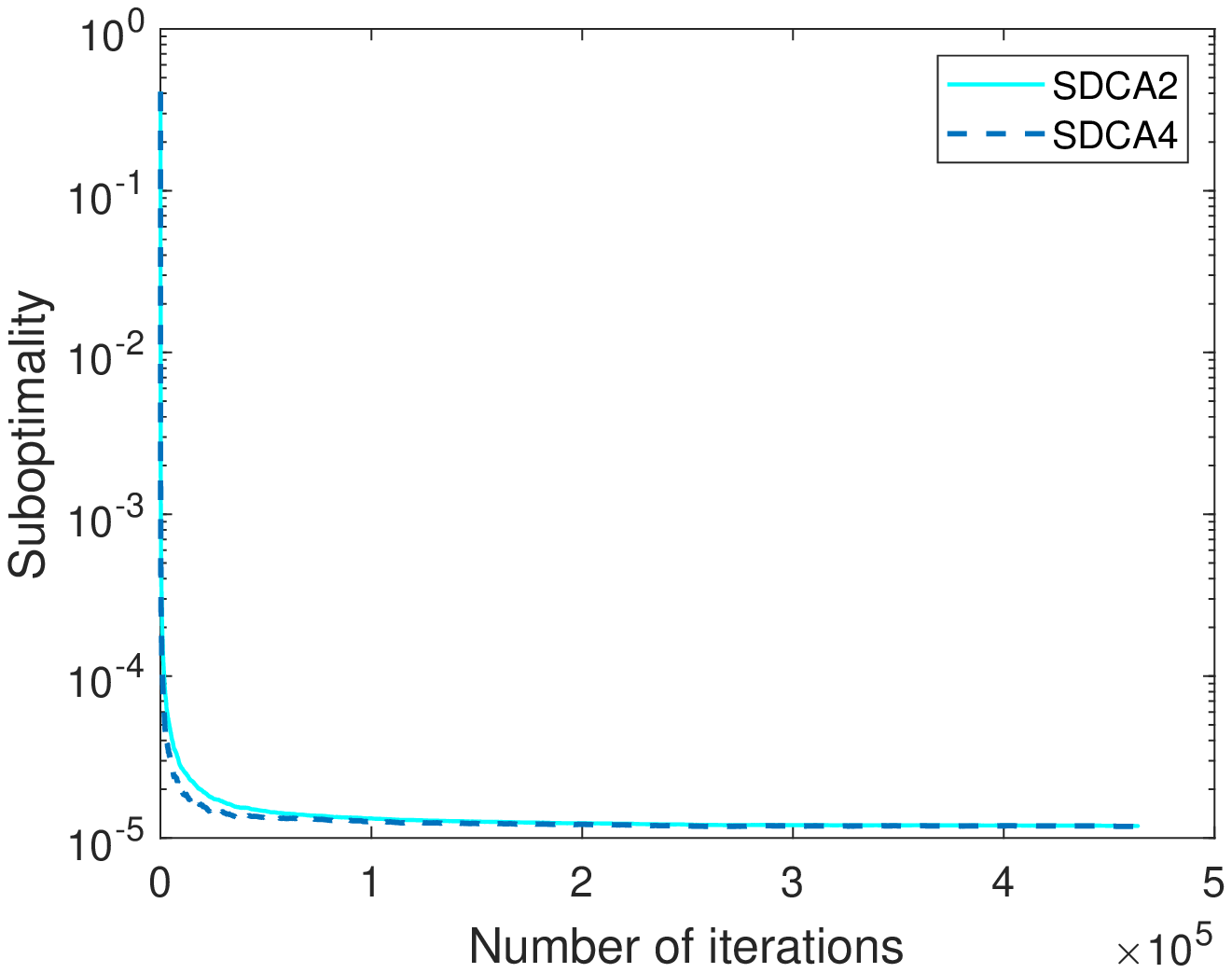} 
	}
	\hspace{-10pt}
	\subfigure[\texttt{SensIT Vehicle}]{	
	\includegraphics[width=0.24\textwidth]{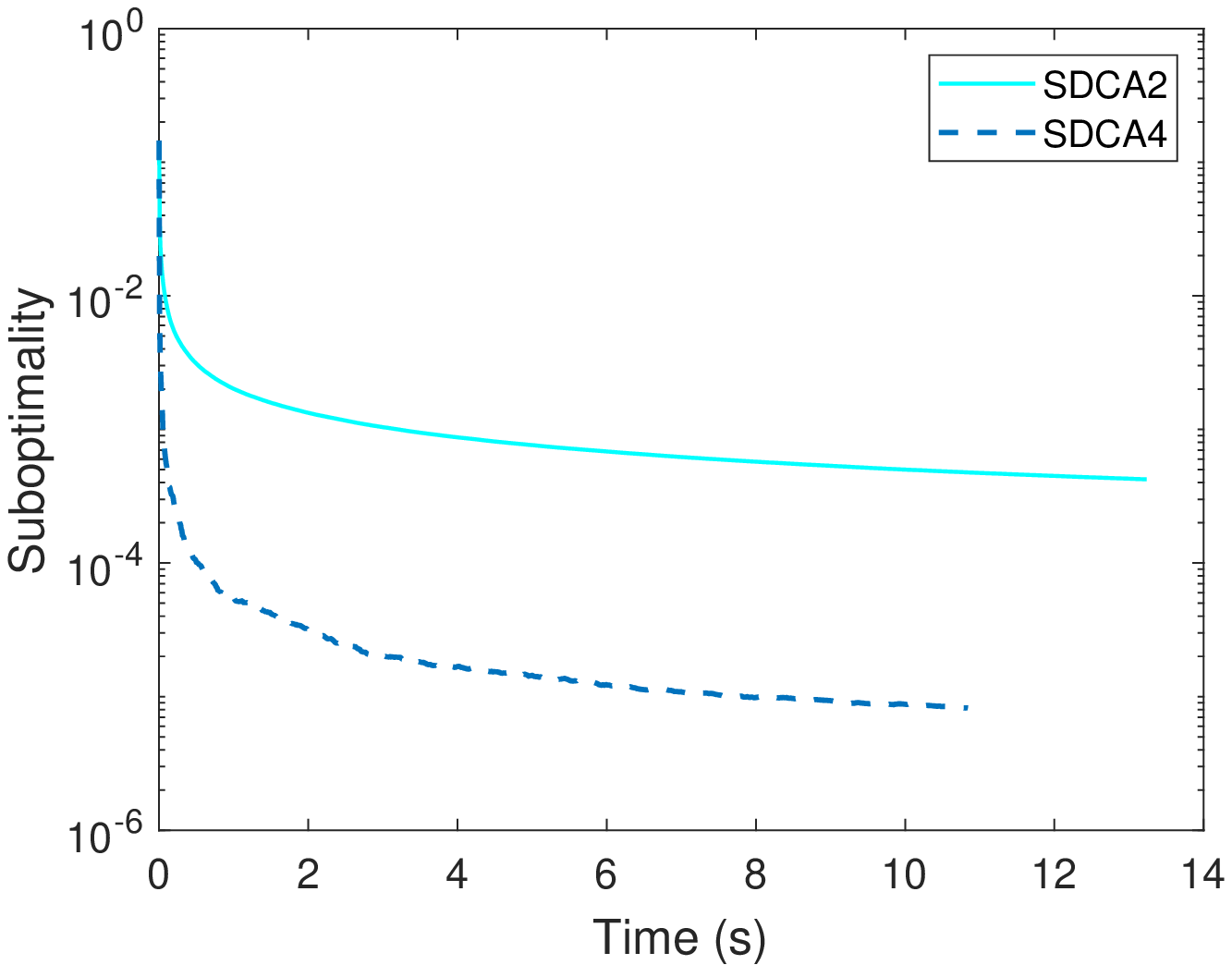} 
	}
	\hspace{-10pt}
	\subfigure[\texttt{shuttle}]{	
	\includegraphics[width=0.24\textwidth]{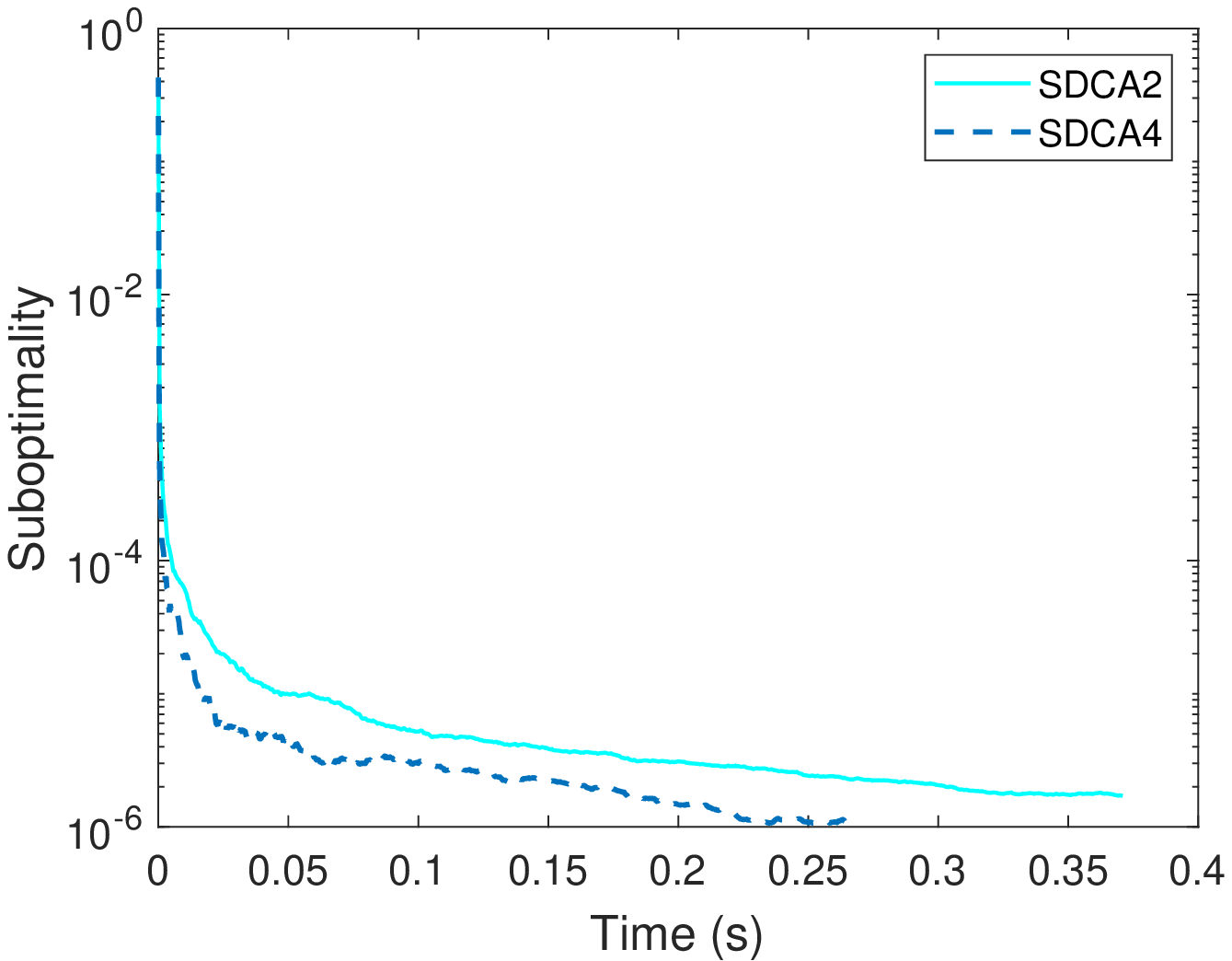} 
	}
	\hspace{-10pt}
	\subfigure[\texttt{protein}]{	
	\includegraphics[width=0.24\textwidth]{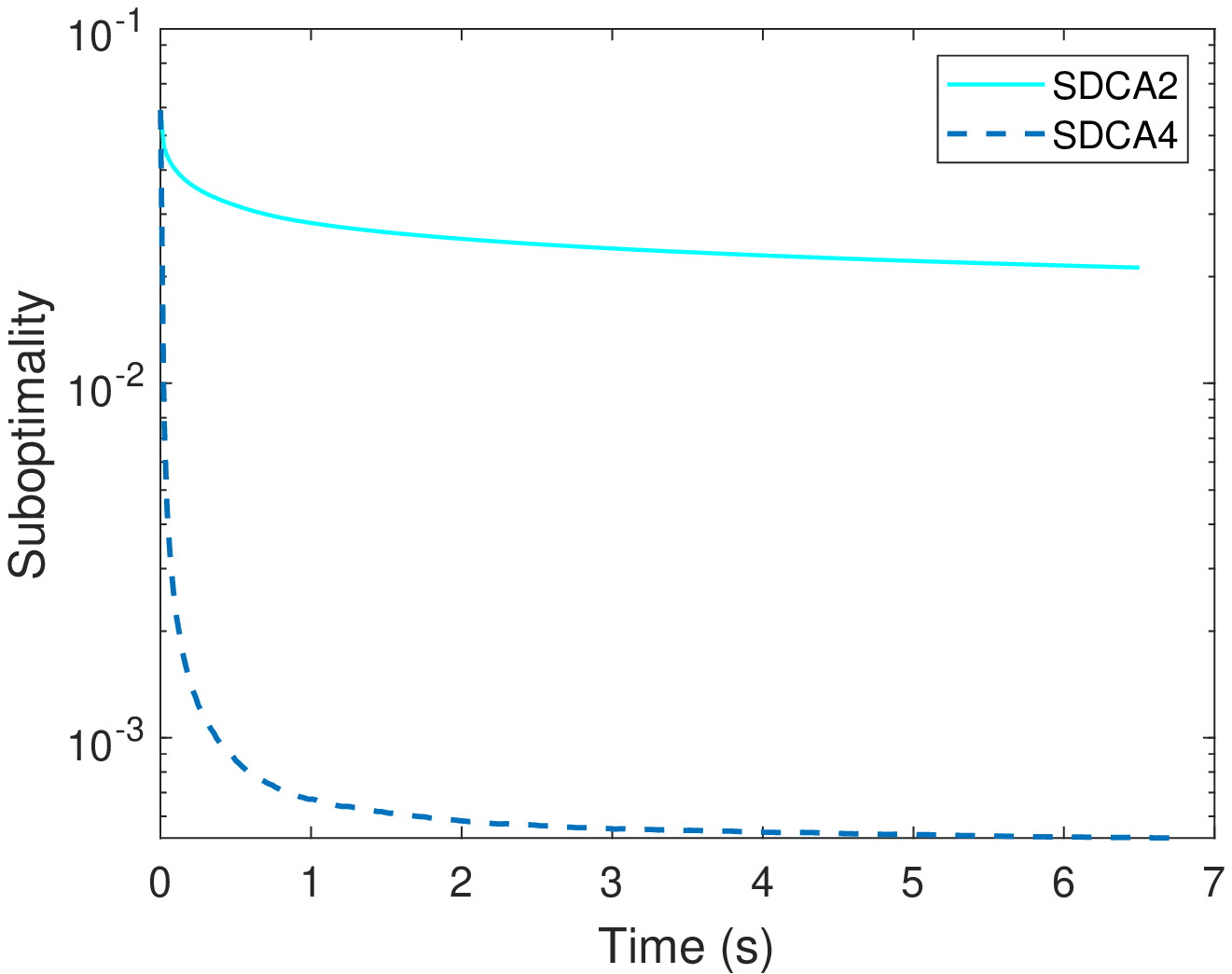} 
	}
	\hspace{-10pt}
	\subfigure[\texttt{YearPredictionMSD}]{	
	\includegraphics[width=0.24\textwidth]{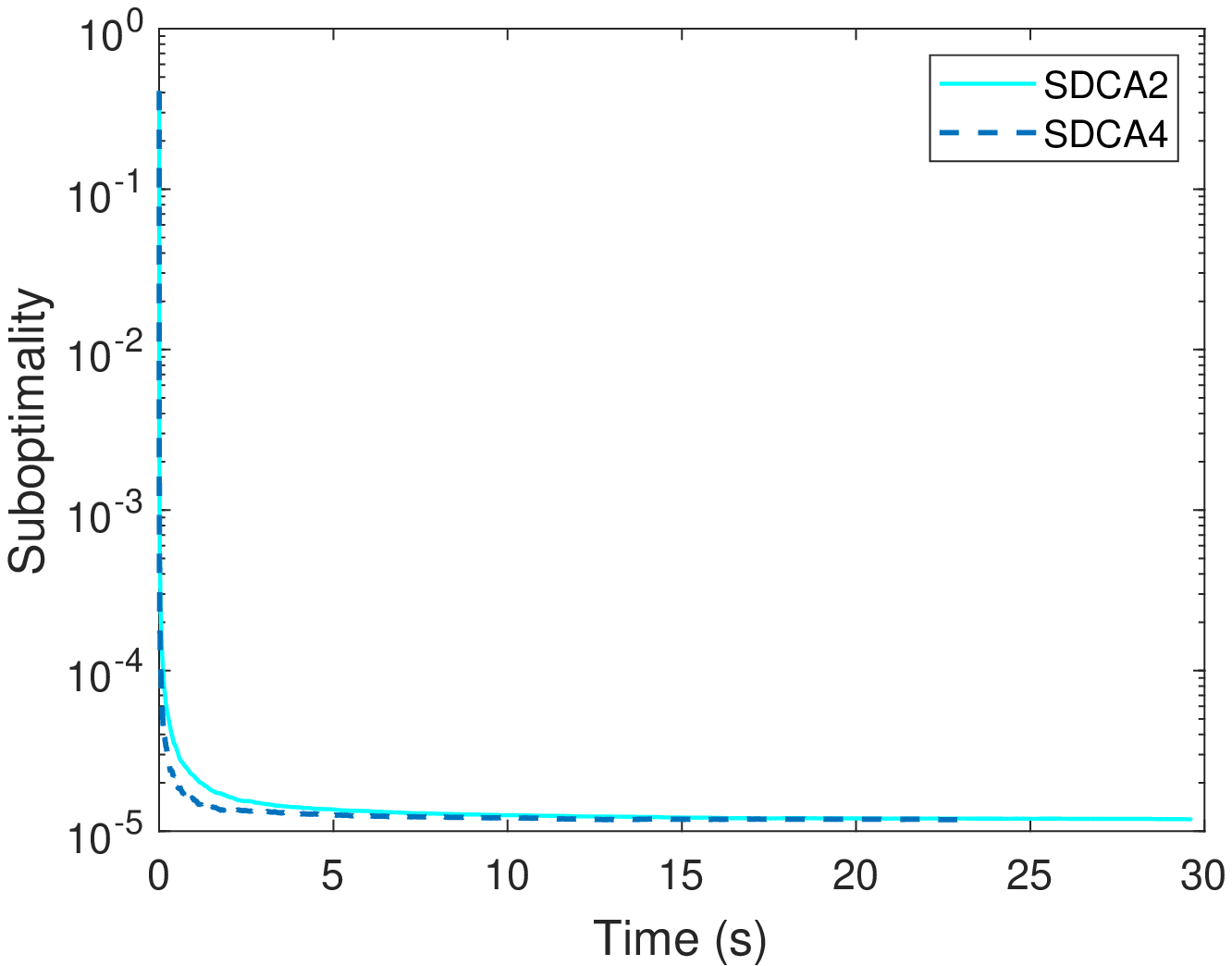} 
	}
     \caption{The performance of SDCA2 and SDCA4 on a ``really" stochastic formulation}
    \label{fig2}
\end{figure}

Though this function is convex, it also has a very natural ``false" DC decomposition with two DC components $\tilde{g}$ and $\tilde{h}$. This observation motivates us to apply deterministic DCA to minimize this function with the stopping criterion being set as $\Vert x^{k+1}-x^k \Vert < 10^{-3}$. Figure \ref{fig2} presents our experimental results. Unlike the experiment 1 where SDCA3,4 does not gain the competitive edge over SDCA1,2, it is clearly observed that SDCA4 outperforms SDCA2, which indicates that the extra effort of recomputing subgradients with respect to past samples actually pays off.

\begin{figure}
     \centering
        \subfigure[\texttt{SensIT Vehicle}]{	
	\includegraphics[width=0.24\textwidth]{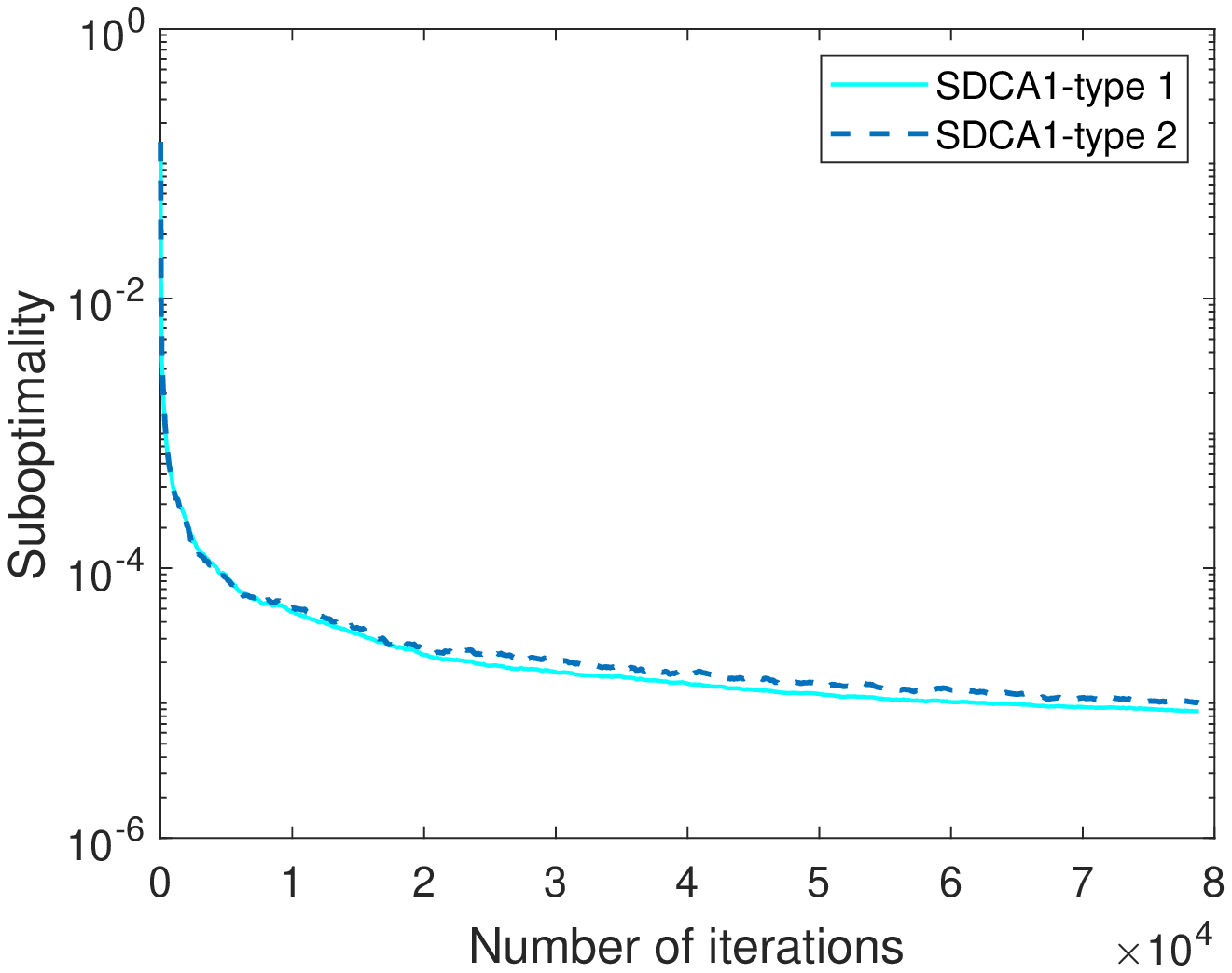} 
	}
	\hspace{-10pt}
	\subfigure[\texttt{shuttle}]{	
	\includegraphics[width=0.24\textwidth]{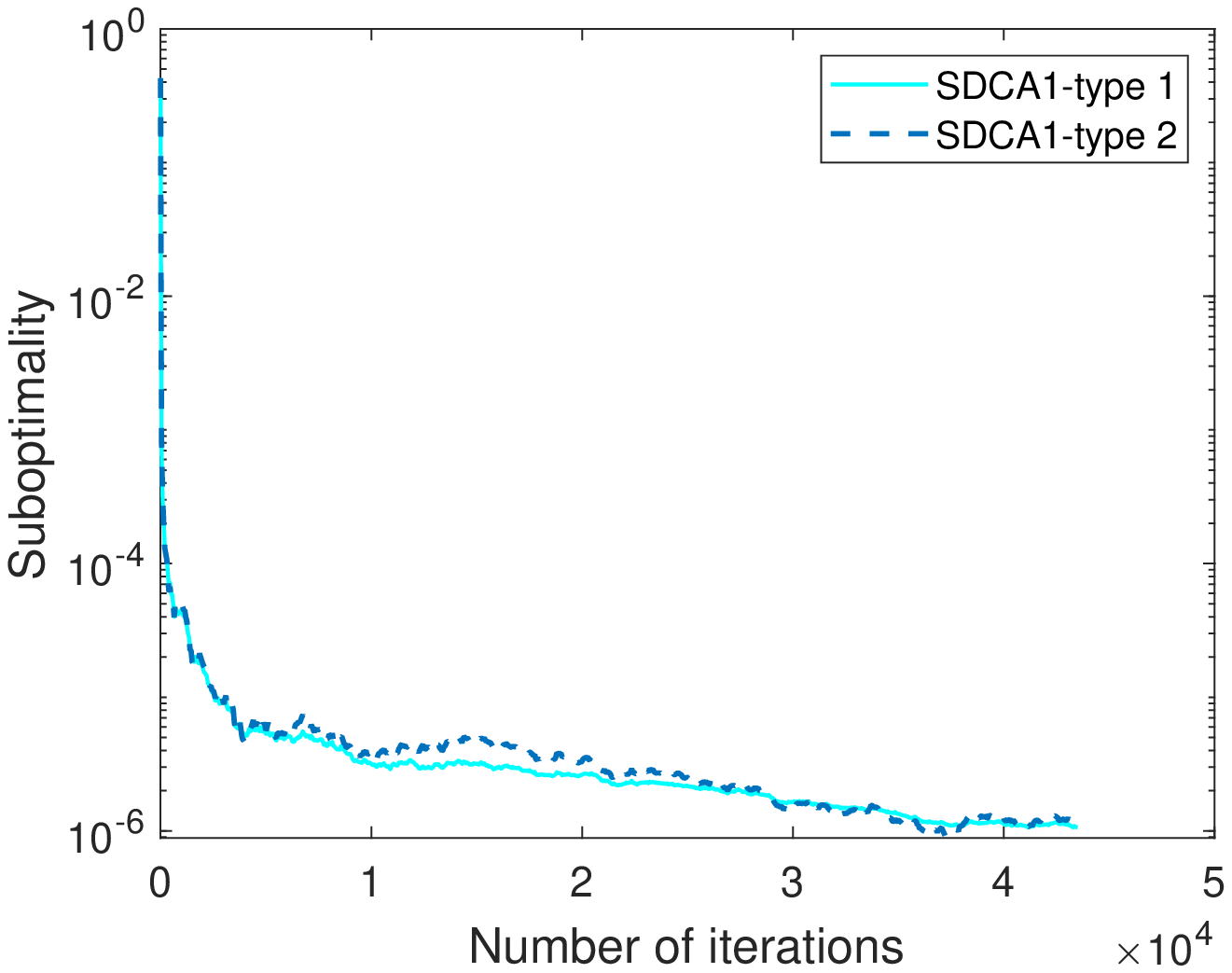} 
	}
	\hspace{-10pt}
	\subfigure[\texttt{protein}]{	
	\includegraphics[width=0.24\textwidth]{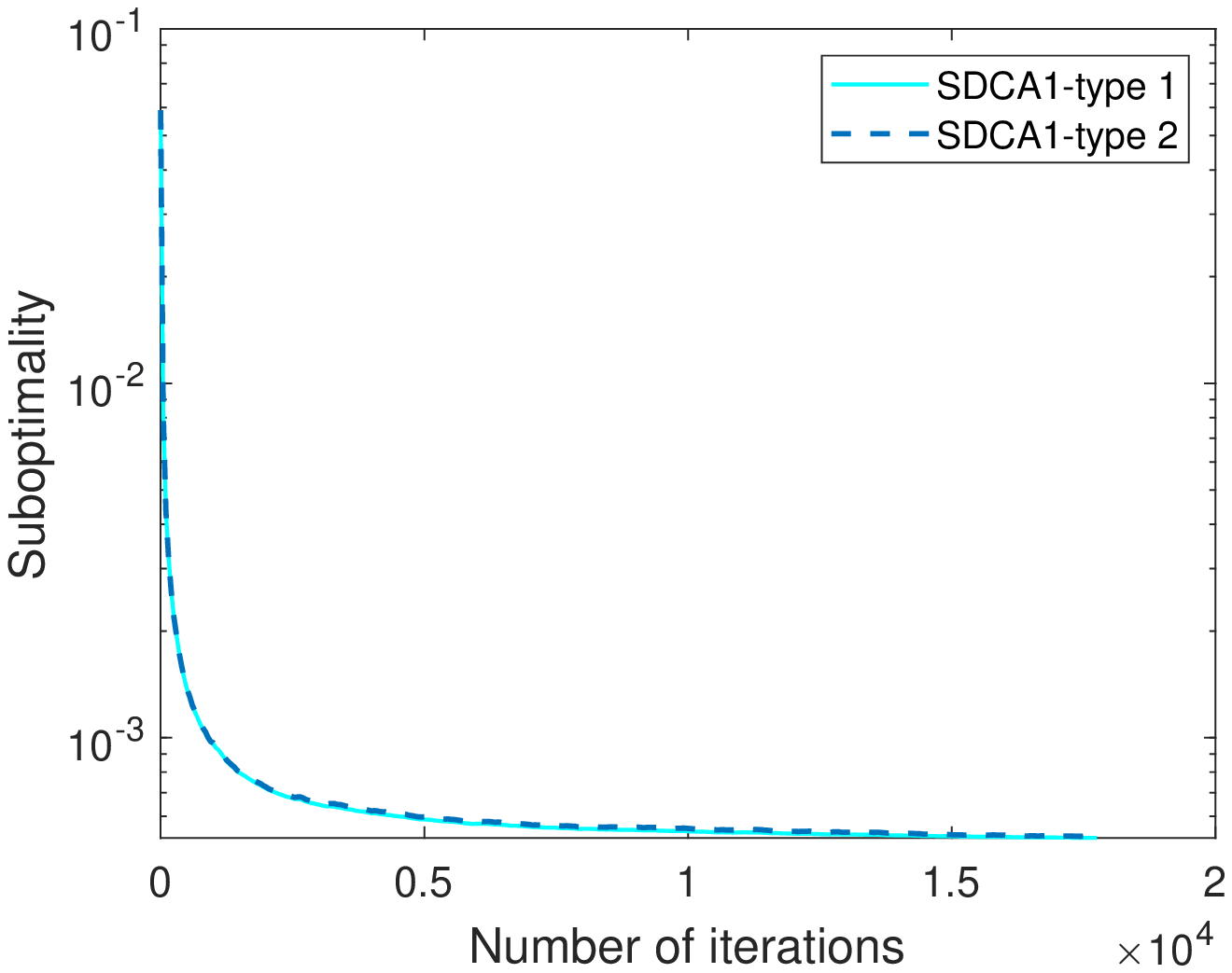} 
	}
	\hspace{-10pt}
	\subfigure[\texttt{YearPredictionMSD}]{	
	\includegraphics[width=0.24\textwidth]{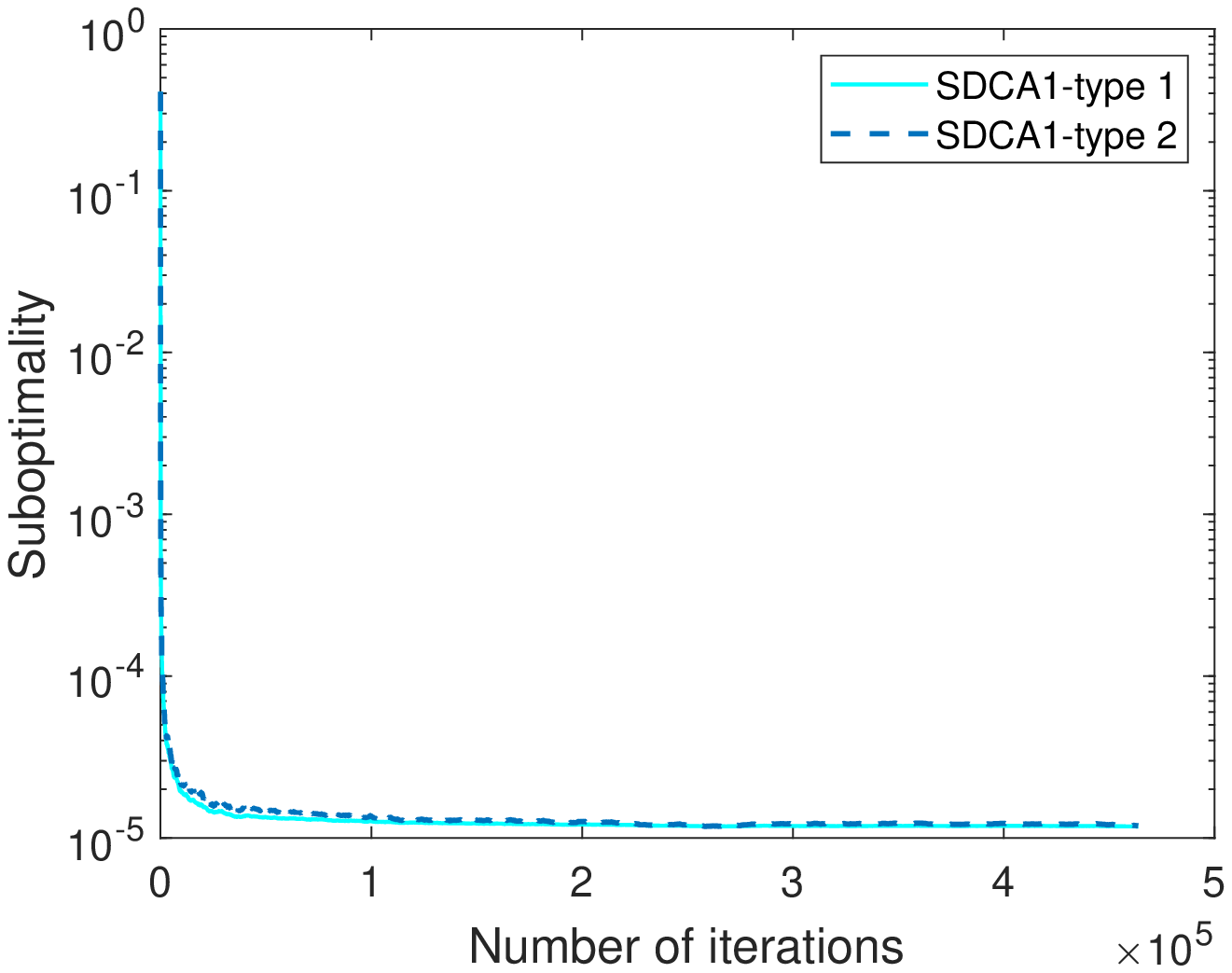} 
	}
	\hspace{-10pt}
	\subfigure[\texttt{SensIT Vehicle}]{	
	\includegraphics[width=0.24\textwidth]{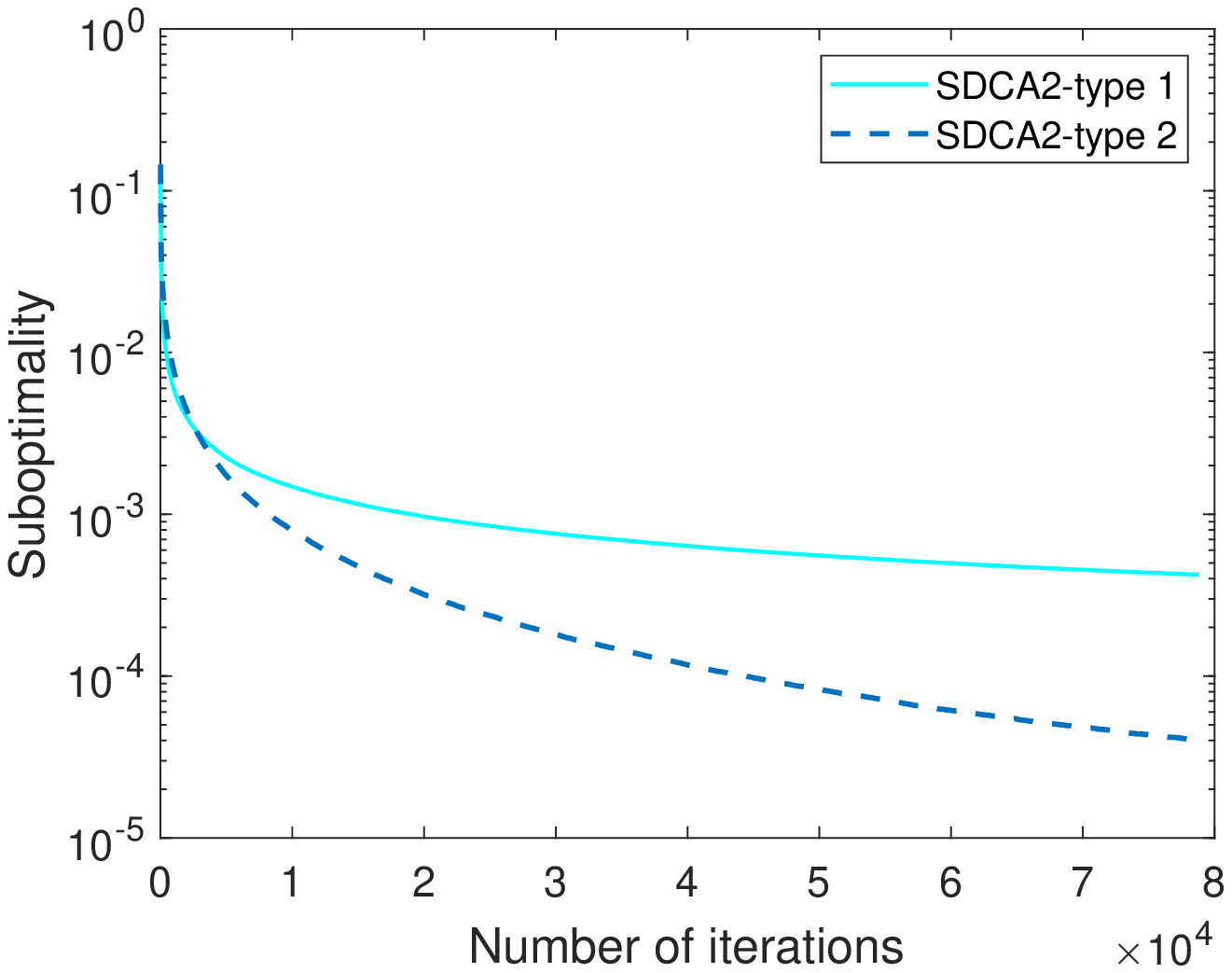} 
	}
	\hspace{-10pt}
	\subfigure[\texttt{shuttle}]{	
	\includegraphics[width=0.24\textwidth]{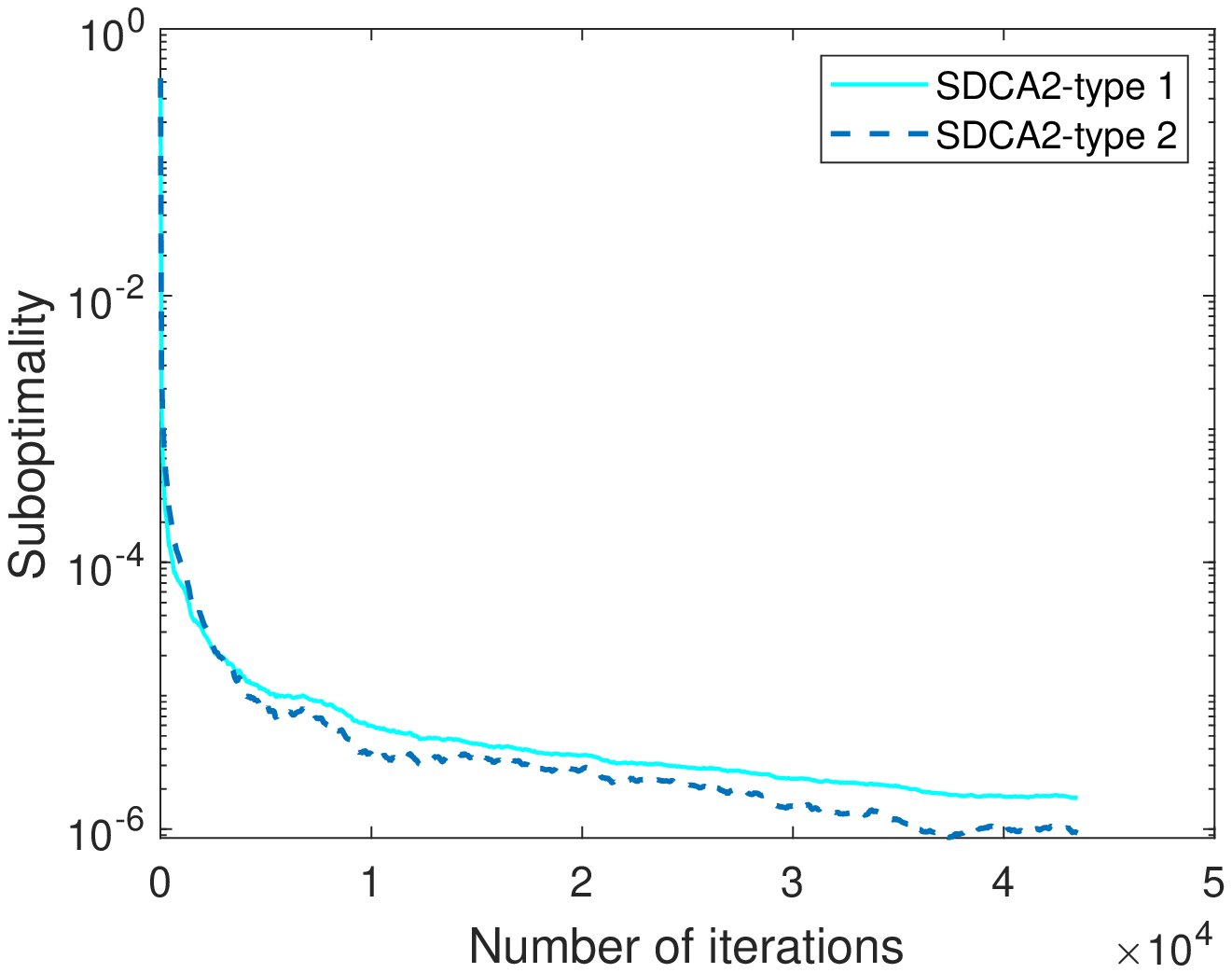} 
	}
	\hspace{-10pt}
	\subfigure[\texttt{protein}]{	
	\includegraphics[width=0.24\textwidth]{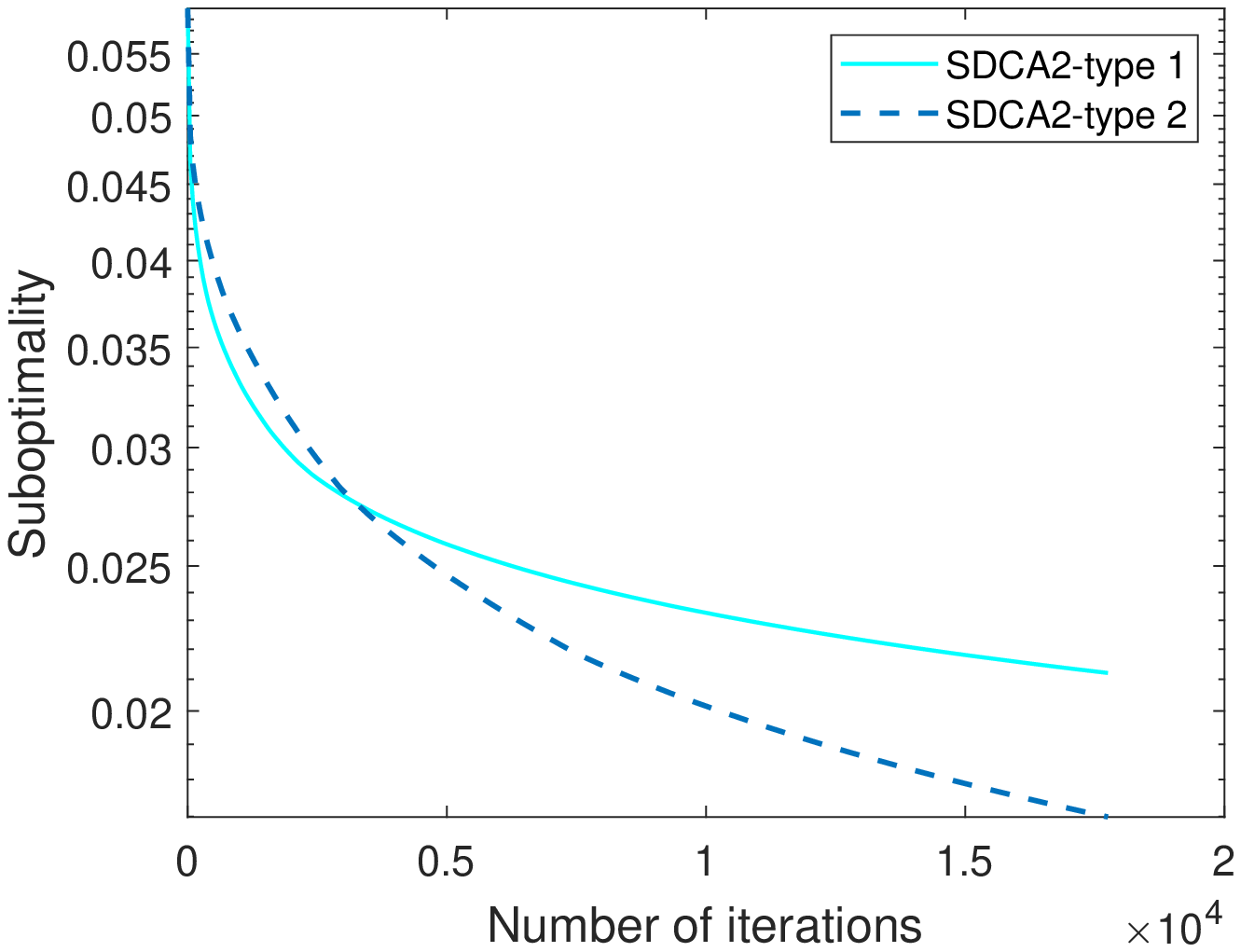} 
	}
	\hspace{-10pt}
	\subfigure[\texttt{YearPredictionMSD}]{	
	\includegraphics[width=0.24\textwidth]{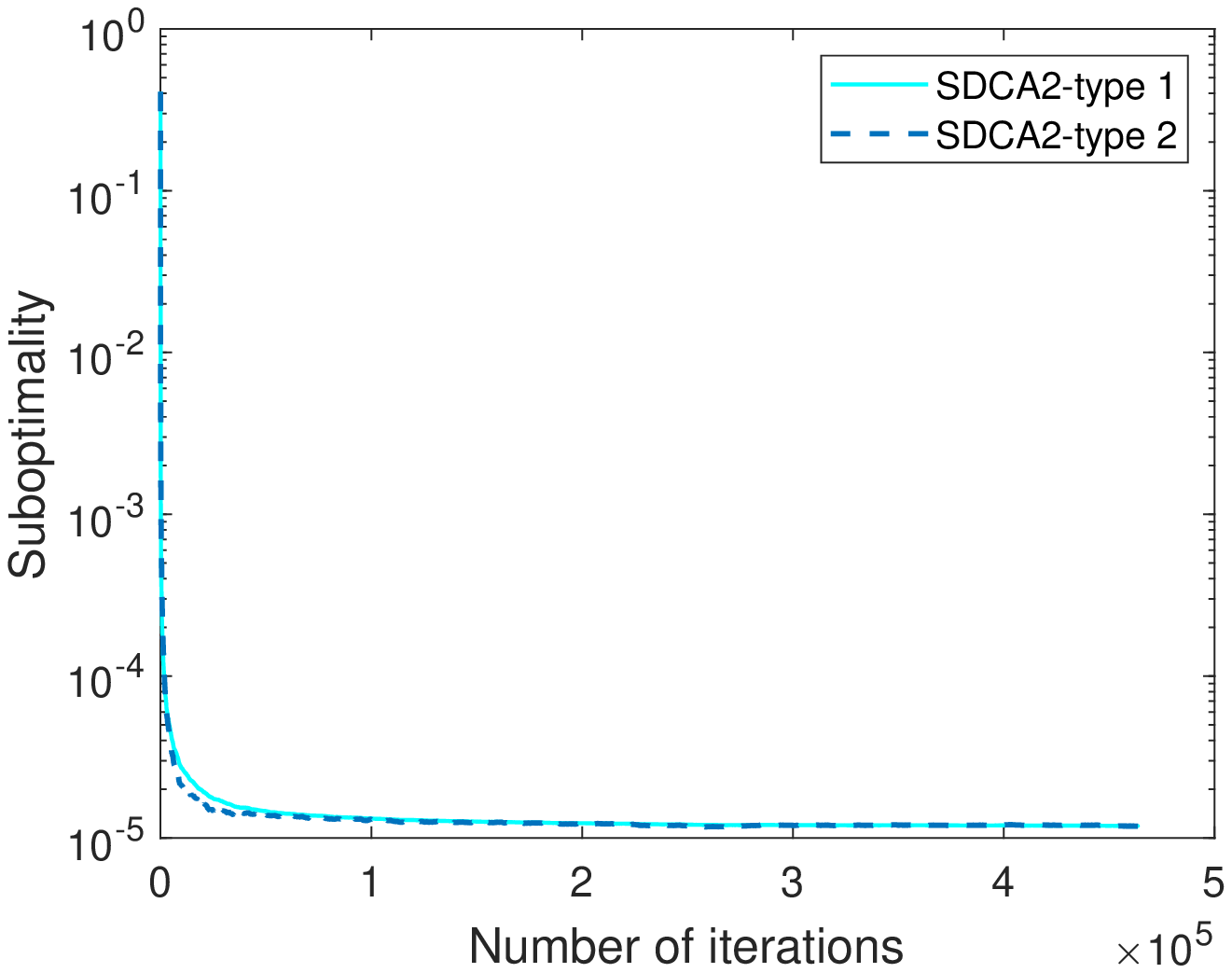} 
	}

     \caption{The effects of weights on the behaviors of SDCA1 and SDCA2}
     \label{fig3}
\end{figure}

In the last experiment, our aim is to observe the effects of the sequence of weights on behaviors of SDCA schemes. To be specific, motivated by \emph{Remark 1}, we want to see if the delicate choice $\alpha_0 =1,  \alpha_k = a^{(k+1)^{\alpha}}-a^{k^{\alpha}}, k \in \mathbb{N}_*$ practically improves the performance of SDCA1 and SDCA2 compared to a more natural choice which is a sequence of equal weights. In the experiment, we set $a = 5, \alpha = 0.2$ which are neutral parameters in line with our theoretical discussion. Figure \ref{fig3} illustrates the performance of SDCA1 and SDCA2 with these two types of weights, where type-1 means equal weights and type-2 indicates the other one. We have an observation that while the performances of SDCA1 with two types of weights are quite similar, the second type of weights really has a positive effect on the behaviors of SDCA2 as it boosts the performance of SDCA2 to obtain better objective values at a faster speed.

\section{Concluding remarks}
In this paper, we have proposed two variants of SDCA (including four algorithms) and have established the convergence results for these algorithms. In Algorithms 1 and 2, the realized samples as well as subgradients from the past iterations are inherited. The convergence rate of Algorithms 3 and 4 in which the subgradients with respect to all past samples  are updated at the current iteration, is considerably faster than Algorithms 1 and 2. We then conducted numerical experiments to study the algorithms' behaviors. Interestingly, the theoretical analysis and numerical performance agree at some points. Moreover, an important procedure in the proposed SDCA is to solve convex optimization subproblems. For  this purpose,   existing stochastic convex optimization approaches, such as the stochastic proximal subgradient methods, could be used. The practical convergence rate of the algorithms depends strictly on the methods dealing with those convex subprograms. The convergence analysis of
the proposed SDCA according to the stochastic methods for solving the convex subproblems 
could be a challenge for further research.  
 \section*{Acknowledgment}
Part of this work has been done during the visit of the second author at LGIPM, University of Lorraine. The second author thanks the University of Lorraine for his financial support and thanks LGIPM for the hospitality.

\section{ Appendix. Proofs of Lemmas \ref{STLN} and \ref{ULLN}}
\textit{Proof of Lemma \ref{STLN}.} The idea of the proof is standard, as the one for the classical strong law of large number (see \cite{Durrett}). We prove the lemma firstly for the case where $\mathbb{E}X_0^4<+\infty$ and $\gamma>1/2.$ Note that $\mathbb{E}X_i^4 = \mathbb{E}X_0^4$ for all $i \in \mathbb{N}$ since $X_0, X_1, \ldots$ are i.i.d. By considering $X_i-\mu$ instead of $X_i$, and as $\mathbb{E}X_i^4<+\infty$ implies $\mathbb{E}(X_i-\mu)^4<+\infty$ (thanks to H\"{o}lder inequality $\mathbb{E}\vert X_i \vert^p \leq (\mathbb{E}X_i^4)^{p/4}$ for $0 \leq p \leq 4$), we can assume that $EX_i=0,$ for all $i\in \N.$ Setting $S_k:=\sum_{i=0}^k\alpha_iX_i,$ since
$$\mathbb{E}S_k^4=\sum_{0\le i,j,l,m\le k}\alpha_i\alpha_j\alpha_l\alpha_m\mathbb{E}(X_iXjX_lX_m),$$
and by independence,
$\mathbb{E}(X_i^3X_j)=\mathbb{E}(X_i^2X_jX_l)=\mathbb{E}(X_iX_jX_lX_m)=0$ for all $0\le i\not=j\le k;$ $1\le i\not=j\not=m$ and all $1\le i\not=j\not= l\not=m\le k,$ one has 
$$\begin{array}{ll}
\mathbb{E}S_k^4&=\mathbb{E}X_0^4\sum_{i=0}^k\alpha_i^4+(\mathbb{E}X_0^2)^2\sum_{0\le i\not=j\le k}\alpha_i^2\alpha_j^2\\
&=\mathbb{E}X_0^4\sum_{i=0}^k\alpha_i^4+(\mathbb{E}X_0^2)^2[(\sum_{i=0}^k\alpha_i^2)^2-\sum_{i=0}^k\alpha_i^4]\le C(\sum_{i=0}^k\alpha_i^2)^2,
\end{array}$$
for some $C>0.$ Therefore, the Chebyshev inequality implies, for any $\varepsilon>0,$
$$\mathbb{P}\left(|S_k|\ge \varepsilon\sum_{i=0}^k\alpha_i\right)\le \varepsilon^{-4}\frac{\mathbb{E}S_k^4}{(\sum_{i=0}^k\alpha_i)^4}\le \frac{C\varepsilon^{-4} (\sum_{i=0}^k\alpha_i^2)^2}{(\sum_{i=0}^k\alpha_i)^4}\le C\varepsilon^{-4} /k^{2\gamma},$$
Consequently,
$\sum_{k=1}^\infty\mathbb{P}\left(|S_k|\ge \varepsilon\sum_{i=0}^k\alpha_i\right)\le C\varepsilon^{-4} \sum_{k=1}^\infty1/k^{2\gamma}<+\infty.$
As $\varepsilon>0$ is arbitrary, by the Borel-Cantelli Lemma, one has $\frac{S_k}{\sum_{i=0}^k\alpha_i}\rightarrow 0\;\; \mbox{a.s.}$
 \vskip 0.1cm
For the second case, by setting $X_k^+=\max\{X_k,0\}$ and $X_k^-=\max\{-X_k,0\},$ then $X_k^+,X_k^-\ge 0;$ $X_k=X_k^+-X_k^-;$ $\mathbb{E}X_k=\mathbb{E}X_k^+-\mathbb{E}X_k^-$, and 
$\mathbb{E}(X^+_k)^2,\mathbb{E}(X^-_k)^2\le \mathbb{E}X_k^2<+\infty.$ We note that $\{X_k^+\}_{k \in \mathbb{N}}$ are i.i.d. random variables; likewise, $\{X_k^-\}_{k \in \mathbb{N}}$ are i.i.d. random variables. So it is enough to prove the lemma for the case $X_k\ge 0.$  Let $A_k=\sum_{i=0}^k\alpha_i$ and $S_k=\sum_{i=0}^k\alpha_iX_i,$ $k=0,1,....$ Then $\mathbb{E}S_k=\sum_{i=0}^k\alpha_i\mathbb{E}X_i=A_k\mu,$ and $\mathbb{D}S_k=\sum_{i=0}^k\alpha_i^2 \mathbb{D}X_0$ (note that i.i.d. random variables have the same mean and the same variance). The Chebyshev inequality implies, for any $\varepsilon>0,$ 
$$\mathbb{P}\left(|S_k-A_k\mu|>\varepsilon A_k\right)\le \varepsilon^{-2}\frac{\mathbb{D}S_k}{A_k^2}=\varepsilon^{-2}\mathbb{D}X_0\frac{\sum_{i=0}^k\alpha_i^2}{A_k^2}\le \varepsilon^{-2}\mathbb{D}X_0\frac{N}{k^\gamma}.$$
Let $l_k=[k^{2/\gamma}],$ $k\in\N$ be the integer part of $k^{2/\gamma}.$ One has
$$\sum_{k=1}^\infty\mathbb{P}\left(|S_{l_k}-A_{l_k}\mu|>\varepsilon A_{l_k}\right)\le \varepsilon^{-2}N\mathbb{D}X_0 \sum_{k=1}^\infty1/[k^{2/\gamma}]^\gamma<+\infty.$$
Thanks to the Borel-Cantelli Lemma, since $\varepsilon>0$ is arbitrary, one obtains $S_{l_k}/A_{l_k}\rightarrow \mu$ almost surely. To show $S_{k}/A_{k}\rightarrow \mu$ a.s.,  for each $k=0,1,..,$ picking $l_{m(k)}$ such that $l_{m(k)}\le k<l_{m(k)+1},$ then 
$\frac{S_{l_{m(k)}}}{A_{l_{m(k)+1}}}\le \frac{S_{k}}{A_{k}}\le\frac{S_{l_{m(k)+1}}}{A_{l_{m(k)}}}. $
Since $\frac{l_{m(k)+1}} {l_{m(k)}}\to 1,$ $\frac{A_{l_{m(k)+1}}}{A_{l_{m(k)}}}\to 1,$ therefore
$\lim_{k\to\infty}\frac{S_{l_{m(k)}}}{A_{l_{m(k)+1}}}=\frac{S_{l_{m(k)+1}}}{A_{l_{m(k)}}}=\mu \quad \mbox{a.s.},$
implying the desired conclusion $S_{k}/A_{k}\rightarrow \mu$ almost surely as $k\to\infty.$\hfill{$\Box$}
\vskip 0.1cm
\noindent\textit{Proof of Lemma \ref{ULLN}.} The proof is similar to the one of Lemma B2 in \cite{EN-13}. By using symmetrization arguments and Rademacher averages as in Lemma 2.3.6 \cite{Wa-We}, one obtains the following estimate, for $k\in\N_*,$
\begin{equation}\label{Symm estim}\mathbb{E}\max_{x\in X}\left|\frac{1}{\sum_{i=0}^k\alpha_i}\sum_{i=0}^k\alpha_i f(x,s^ i)-\mathbb{E}_sf(x,s)\right|\le 2 \mathbb{E}R_k(f,\alpha, X),
\end{equation}
where, 
$R_k(f,\alpha,X):=\mathbb{E}_\sigma\sup_{x\in X}\frac{1}{\sum_{i=0}^k\alpha_i}\left|\sum_{i=0}^k\sigma_i\alpha_if(x,s^ i)\right|;$
$\{\sigma_i\}$ are i.i.d. random variables taking values $\pm 1$ with probability $1/2$ each. For a set $A\subseteq\R^{k+1},$ denote by $R_k(\alpha,A)$ the Rademacher average of $A$ with respect to $\{\alpha_i\}$:
$$R_k(\alpha,A):=\mathbb{E}_\sigma\sup_{a\in A}\frac{1}{\sum_{i=0}^k\alpha_i}\left|\sum_{i=0}^k\sigma_i\alpha_ia_i\right|,$$ 
where $a=(a_0,a_1,...,a_k)\in\R^{k+1}.$
The following lemma gives an upper bound  of $R_k(\alpha,A),$ which generalizes the one in Theorem 3.3 in \cite{BBL}.
\begin{lemma}\label{Rade estim} Let $A\subseteq\R^{k+1}$ be a finite set of $N$ elements. One has
$$R_k(\alpha,A)\le \max_{a\in A}\max_{i=0,...,k}|a_i|\frac{\sqrt{2\ln(2N)\sum_{i=0}^k\alpha_i^2}}{\sum_{i=0}^k\alpha_i}.$$
\end{lemma}
\textit{Proof.} Set $A_k=\sum_{i=0}^k\alpha_i.$ By using the Hoeffding inequality, stating that for a zero-mean random variable with values in $[t_1,t_2],$ $\mathbb{E}e^X\le e^{(t_2-t_1)^2/8},$ for any $s>0,$ one has
$$\mathbb{E}e^{\frac{s}{A_k}\sum_{i=0}^k\sigma_i\alpha_ia_i}=\Pi_{i=0}^k\mathbb{E}e^{\frac{s}{A_k}\sigma_i\alpha_ia_i}\le\Pi_{i=0}^ke^{\frac{s^2\alpha^2_ia^2_i}{2A^2_k}}=e^{\frac{s^2}{2A^2_k}\sum_{i=0}^k\alpha^2_ia^2_i},$$
where the first equality is by independence.
Hence, by the Jensen inequality,
$$\begin{array}{ll}
&e^{s\mathbb{E}\max_{a\in A}(1/A_k)\sum_{i=0}^k\sigma_i\alpha_ia_i}\le \mathbb{E}e^{s\max_{a\in A}(1/A_k)\sum_{i=0}^k\sigma_i\alpha_ia_i}\\
&\le \sum_{a\in A}\mathbb{E}e^{\frac{s}{A_k}\sum_{i=0}^k\sigma_i\alpha_ia_i}\le N\max_{a\in A}e^{\frac{s^2}{2A^2_k}\sum_{i=0}^k\alpha^2_ia^2_i},
\end{array}$$
which implies
$\mathbb{E}\max_{a\in A}(1/A_k)\sum_{i=0}^k\sigma_i\alpha_ia_i\le \frac{\ln N}{s}+\max_{a\in A}\frac{s}{2A^2_k}\sum_{i=0}^k\alpha^2_ia^2_i.$
By noting that
$R_k(\alpha,A)=\mathbb{E}\max_{a\in A\cup(-A)}(1/A_k)\sum_{i=0}^k\sigma_i\alpha_ia_i,$
the preceding inequality yields, by considering the set $A\cup(-A)$ instead of $A,$
$$R_k(\alpha,A)\le\frac{\ln (2N)}{s}+\max_{a\in A}\frac{s}{2A^2_k}\sum_{i=0}^k\alpha^2_ia^2_i,$$ 
and by setting $s=A_k\sqrt{2\ln(2N)/\max_{a\in A}\sum_{i=0}^k\alpha^2_ia^2_i},$ it implies the desired estimate:
$$R_k(\alpha,A)\le \max_{a\in A}\frac{\sqrt{2\ln(2N)\sum_{i=0}^k\alpha^2_ia^2_i}}{A_k}\le\max_{a\in A}\max_{i=0,...,k}|a_i|\frac{\sqrt{2\ln(2N)\sum_{i=0}^k\alpha_i^2}}{\sum_{i=0}^k\alpha_i}.$$
\hfill{$\Box$}
\vskip 0.1cm
To end the proof of Lemma \ref{ULLN}, in view of estimate (\ref{Symm estim}), we shall show that for some $c>0,$ for all $k\in \N_*,$
\begin{equation}\label{Rademacher upper bound}
R_k(f,\alpha,X)\le \frac{c(1+\sqrt{\ln \beta_k})}{\beta_k}, 
\end{equation}
where $\beta_k:=\frac{\sum_{i=0}^k\alpha_i}{\left(\sum_{i=0}^k\alpha_i^2\right)^{1/2}}.$
Indeed, assume that $X\subseteq B_R\subseteq\R^n$, a Euclidean ball centered at $0$ with radius $R>0.$ For any $\varepsilon>0,$ there are $N:=N(\varepsilon)\le 4Re^n/\varepsilon$ balls with radius $\varepsilon:$ $B(y^i,\varepsilon),\;y^i\in X,\; i=1,...,N$ covering $X,$ that is, $X\subseteq \bigcup_{i=1}^N B({y^i},\varepsilon)$ (see, e.g., \cite{CS}). Let $Y_\varepsilon:=\{y^1,...,y^N\},$ since $f(\cdot,s)$ are H\"{o}lder continuous on $X$ with constants $L,\gamma>0$, one has 
$\left|R_k(f,\alpha,X)-R_k(f,\alpha,Y_\varepsilon)\right|\le L\varepsilon^\gamma.$
Let $A:=\left\{a^i=(f(y^i,s_0),...,f(y^i,s^k)):\;\;i=0,...,N\right\},$ we see that $R_k(f,\alpha,Y_\varepsilon)=R(\alpha,A),$ therefore Lemma \ref{Rade estim} implies
\begin{align*}
R_k(f,\alpha,Y_\varepsilon)=R(\alpha,A) &\leq \max_{a\in A}\max_{i=0,...,k}|a_i|\frac{\sqrt{2\ln(2N)\sum_{i=0}^k\alpha_i^2}}{\sum_{i=0}^k\alpha_i}\\
&\le M\frac{\sqrt{2(n+\ln(8R/\varepsilon))\sum_{i=0}^k\alpha_i^2}}{\sum_{i=0}^k\alpha_i}
\end{align*}
since $|f(x,s)|\le M$ for all $(x,s)\in X\times \Omega$.

This inequality together with the preceding inequality imply 
$$R_k(f,\alpha,X)\le  L\varepsilon^\gamma + M\frac{\sqrt{2(n+\ln(8R/\varepsilon))\sum_{i=0}^k\alpha_i^2}}{\sum_{i=0}^k\alpha_i},$$
and by taking 
$\varepsilon=\left(\frac{\sqrt{\sum_{i=0}^k\alpha_i^2}}{\sum_{i=0}^k\alpha_i}\right)^{1/\gamma},$ we derive the desired estimate.\hfill{$\Box$}
\end{document}